\documentclass[10pt,twoside]{article}
\usepackage{amsmath}
\usepackage[utf8]{inputenc}
\usepackage{url}
\usepackage{dsfont}
\usepackage{graphicx}
\usepackage{yhmath}
\usepackage[table]{xcolor}
\usepackage{mathrsfs} % permet d'avoir le D cursive
\usepackage{amssymb}
\usepackage{url}
\usepackage{makecell}
\usepackage{arydshln}
\usepackage{multirow}
\usepackage{arydshln}
\usepackage{mathtools}
\usepackage{stmaryrd}  % pour intervalle entiers
 \usepackage{amsthm,amssymb}
\usepackage{hyperref}

\input epsf
\def\limiten{\renewcommand{\arraystretch}{0.5}
\begin{array}[t]{c}\stackrel{}{\longrightarrow} \\
{\scriptstyle n\rightarrow
\infty}\end{array}\renewcommand{\arraystretch}{1}}

\numberwithin{equation}{section}

\newtheorem{thm}{Theorem}[section]

\newtheorem{Corol}[thm]{Corollary}
\newtheorem{Def}[thm]{Definition\rm}

\newtheorem{lem}[thm]{Lemma}

\newcommand{\E}{\ensuremath{\mathbb{E}}}
\newcommand{\R}{\ensuremath{\mathbb{R}}}
\newcommand{\Z}{\ensuremath{\mathbb{Z}}}

\newcommand{\N}{\ensuremath{\mathbb{N}}}

\newcommand{\var}{\ensuremath{\mathrm{Var}}}

\definecolor{grisclair}{gray}{0.9}
\font\dsrom=dsrom10 scaled 1200
\def \ind{\textrm{\dsrom{1}}}
\DeclareMathOperator*{\argmin}{argmin}

 \newcommand{\mk}{ { \mathcal{K}} }
\newcommand{\mx}{ { \mathcal{X}} }

\newcommand{\mz}{ { \mathcal{Z}} }

\renewcommand{\arraystretch}{.8}

\begin{document}
%\date{ }

\title{\bf A general framework for deep learning}
 \maketitle \vspace{-1.0cm}

\begin{center}
      William Kengne$^1$
     % \footnote{Developed within the ANR BREAKRISK: ANR-17-CE26-0001-01 and the  CY Initiative of Excellence (grant "Investissements d'Avenir" ANR-16-IDEX-0008), Project "EcoDep" PSI-AAP2020-0000000013 } 
   and 
     Modou Wade$^2$
     %\footnote{Supported by the MME-DII center of excellence (ANR-11-LABEX-0023-01)} 
 \end{center}

  \begin{center}
  { \it 
  $^1$Université Jean Monnet, ICJ UMR5208, CNRS, Ecole Centrale de Lyon, INSA Lyon, Universite Claude Bernard Lyon 1, Saint-Étienne, France\\  
 $^2$THEMA, CY Cergy Paris Université, 33 Boulevard du Port, 95011 Cergy-Pontoise Cedex, France\\
  E-mail:   william.kengne@univ-st-etienne.fr  ; modou.wade@cyu.fr\\
  }
\end{center}

 \pagestyle{myheadings}
%\markboth{}{ Kengne and Wade}

\markboth{ A general framework for deep learning}{Kengne and Wade}

\medskip

\textbf{Abstract}:
This paper develops a general approach for deep learning for a setting that includes nonparametric regression and classification.
We perform a framework from data that fulfills a generalized Bernstein-type inequality, including independent, $\phi$-mixing, strongly mixing and $\mathcal{C}$-mixing observations.
Two estimators are proposed: a non-penalized deep neural network estimator (NPDNN) and a sparse-penalized deep neural network estimator (SPDNN).
For each of these estimators, bounds of the expected excess risk on the class of H\"older smooth functions and composition H\"older functions are established.
Applications to independent data, as well as to $\phi$-mixing, strongly mixing, $\mathcal{C}$-mixing processes are considered. 
For each of these examples, the upper bounds of the expected excess risk of the proposed NPDNN and SPDNN predictors are derived.
It is shown that both the NPDNN and SPDNN estimators are minimax optimal (up to a logarithmic factor) in many classical settings.

\medskip
 

\medskip

\section{Introduction}
In the past few decades, the study of theoretical guarantees of deep neural networks (DNN) algorithms has received significant attention in the literature; see, for instance \cite{bauer2019deep}, \cite{schmidt2020nonparametric}, \cite{padilla2022quantile}, \cite{jiao2023deep} and the references therein.
Several works established that the DNN estimators can achieve an optimal convergence rate for regression and classification tasks; see for instance \cite{schmidt2020nonparametric}, \cite{kim2021fast}, \cite{ohn2022nonconvex}, \cite{jiao2023deep}, \cite{zhang2024classification}, \cite{fan2024noise} for some results based on independent and identically distributed (i.i.d.) data and \cite{kohler2023rate}, \cite{alquier2025minimax}, \cite{kurisu2025adaptive}, \cite{kengne2025deep} for some recent advances based on dependent observations.
In general, and whatever one deals with independent data or a specific structure of dependence, the convergence rates mainly depend on the concentration inequality that can satisfy the data.
A natural question then arises: given a generalized Bernstein-type inequality satisfies by the observations, can the convergence rate of the DNN estimator be directly deduced?
The answer to this question will be the main thread of this contribution.  

\medskip

Let us consider the training sample $ D_n \coloneqq  \{ Z_1 = (X_1, Y_1), \cdots, Z_n= (X_n, Y_n) \} $ which is a trajectory of a stationary and ergodic process $\{Z_t =(X_t, Y_t), t \in \Z \} $, which takes values in $\mathcal{Z} = \mathcal{X}  \times \mathcal{Y} \subset \R^d \times \R $ (with $d \in \N$), where $\mathcal{X}$ is the input space and $\mathcal{Y}$ the output space. 
In the sequel, we deal with a broad class of loss functions $\ell: \R \times  \mathcal{Y} \to [0, \infty)$.
Let $\mathcal{F}$ denote the set of measurable functions from $\mx$ to $ \mathcal{Y}$. 
For any predictor $ h\in \mathcal{F}$, its risk is defined as follows:
\begin{equation}\label{equa_def_risk}
 R(h) = \E_{Z_0} [\ell (h(X_0), Y_0)], ~ \text{with}~ Z_0 = (X_0, Y_0).  
\end{equation}
The target predictor (assumed to exist) $h^{*} \in \mathcal{F} $ satisfies:
\begin{equation}\label{best_pred_F}
R(h^{*}) =  \underset{h \in \mathcal{F}}{\inf}R(h).
\end{equation}
The excess risk of a predictor $h \in \mathcal{F}$, is given by:
\begin{equation}\label{equa_excess_risk}
\mathcal{E}_{Z_0}(h) = R(h) - R(h^{*}), ~ \text{with} ~ Z_0 = (X_0, Y_0).
\end{equation} 

\medskip

Let $\sigma: \R \to \R$ be an activation function. We deal with the class of DNN predictors with network architecture $ (L, \textbf{p}) $, where $L \in \N$ denotes the number of hidden layers (or depth) and $\mathbf{p} = (p_0, p_1,\cdots, p_{L+1}) \in \N^{L+2}$ is the width vector.
Such functions are of the form:
\begin{equation}\label{h_equ1}
h: \R^{p_0} \rightarrow \R^{p_{L+1}}, \;   x\mapsto h(x) = A_{L+1} \circ \sigma_{L} \circ A_{L} \circ \sigma_{L-1} \circ \cdots \circ \sigma_1 \circ A_1 (x),
\end{equation} 
where for any $j, \cdots, L+1$, $ A_j: \R^ {p_ {j -1}} \to \R^ {p_j} $ is a linear affine map, defined by $A_j (x) \coloneqq W_j x + \textbf{b}_j$,  for given  $p_ {j}\times p_{j-1}$  weight matrix   $ W_j$   and a shift vector $ \textbf{b}_j \in \R^ {p_j} $, and $\sigma_j: \R^{p_j} \rightarrow \R^ {p_j} $ is a nonlinear element-wise activation map, defined for all $z=(z_1, \cdots, z_{p_j})^T$ by $\sigma_j (z) = (\sigma(z_1), \cdots, \sigma(z_{p_j}))^{T} $, and $^T$ denotes the transpose.   
Let us denote by,
\begin{equation} \label{def_theta_h}
\theta(h) \coloneqq \left(vec(W_1)^ {T}, \textbf{b}^{T}_{1}, \cdots,  vec(W_{L + 1})^{T} , \textbf{b}^ {T}_{L+1} \right)^{T}, 
 \end{equation} 
 the vector of parameters for a DNN function of the form  (\ref{h_equ1}), where $ vec(W)$ is the vector obtained by concatenating the column vectors of the matrix $W$. 
 In the following, we consider an activation function $ \sigma: \R \to \R$ and let us denote by $\mathcal{H}_{\sigma, p_0, p_{L+1}} $, the class of DNN predictors with $p_0$-dimensional input and $p_{L+1} $-dimensional output.
For the framework considered here, $p_0 = d$ and $p_ {L + 1} = 1$.
For a DNN $h$ as in (\ref{h_equ1}), denote by depth($h$)$=L$ and width($h$) = $\underset{1\leq j \leq L} {\max} p_j $ its depth and width respectively.  For any positive constants $L, N, B, F, S > 0$, set
\begin{equation}\label{no_constrained_dnn_class}\mathcal{H}_{\sigma}(L, N, B, F)  \coloneqq \big\{ h: \mathcal{H}_{\sigma, d, 1}: \text{depth}(h)\leq L, \text{width}(h)\leq N, \|\theta(h)\|_{\infty} \leq B, \| h \|_{\infty, \mathcal{X}}  \leq F \big\}, 
\end{equation}
 and
\begin{equation}\label{constrained_dnn_class}
\mathcal{H}_{\sigma}(L, N, B, F, S)  \coloneqq \left\{h\in \mathcal{H}_{\sigma}(L, N, B, F) \; :  \; \| \theta(h) \|_0 \leq S 
  \right\}, 
\end{equation}
where $\| x \|_0 = \sum_{i=1}^p \ind(x_i \neq 0), ~ \| x\|_{\infty} = \underset{1 \leq i  \leq p}{\max} |x_i |$ for all $x=(x_1,\ldots,x_p)^{T} \in \R^p$ ($p \in \N$) 
and $\| h \|_{\infty, \mathcal{X}} = \sup_{x \in \mx} |h(x)|$ for all $ h \in \mathcal{F}$.
$\mathcal{H}_{\sigma}(L, N, B, F, S)$ is a class of sparsity constrained DNN with sparsity level $S > 0$.

\medskip

%The goal is to derive bounds of the expected excess risk of the DNN estimators considered and to provide their convergence rates.

\medskip

We focus on a class of sparsity-constrained DNN $\mathcal{H}_{\sigma}(L_n, N_n, B_n, F_n, S_n)$ defined in (\ref{constrained_dnn_class}), with an activation function $\sigma$ and a suitably chosen network architecture $(L_n, N_n, B_n, F_n, S_n)$. 
The  empirical minimizer over the class of DNN functions $\mathcal{H}_{\sigma}(L_n, N_n, B_n, F_n, S_n)$, also called non-penalized DNN (NPDNN) estimator is defined by: 
\begin{equation}\label{NP_DNNs_Estimators}
\widehat{h}_{n, NP} = \underset{h \in  \mathcal{H}_{\sigma}(L_n, N_n, B_n, F_n, S_n)}{\argmin} \left[ \dfrac{1}{n} \sum_{i=1}^{n} \ell(h(X_i), Y_i)\right].
\end{equation}
%
%
%We focus on a DNN class $\mathcal{H}_{\sigma}(L_n, N_n, B_n, F_n)$ defined in (\ref{no_constrained_dnn_class}), with an activation function $\sigma$ and a suitably chosen network architecture $(L_n, N_n, B_n, F_n)$. 
% 
%  
%The goal is to build, from the training sample $D_n$, a DNN predictor $\widehat{h}_n \in \mathcal{H}_{\sigma}(L_n, N_n, B_n, F)$ with suitably chosen network architecture $(L_n, N_n, B_n, F)$, such that, its expected excess risk tends to zero as the sample size tends to infinity.
%
%In this contribution, we consider the sparse-penalized regularization.
%
We also perform a sparse-penalized DNN (SPDNN) estimator defined by: 
\begin{equation}\label{sparse_DNNs_Estimators}
\widehat{h}_{n, SP} = \underset{h \in  \mathcal{H}_{\sigma}(L_n, N_n, B_n, F_n)}{\argmin} \left[ \dfrac{1}{n} \sum_{i=1}^{n} \ell(h(X_i), Y_i) + J_n(h) \right],
\end{equation}
where for any neural network estimator $h \in \mathcal{H}_{\sigma}(L_n, N_n, B_n, F)$, the sparse penalty term $J_{n}(h) $ fulfills: 
\begin{equation}\label{equa_penalty_term_v1}
J_n(h) = \sum_{j=1}^{p} \pi_{\lambda_n, \tau_n} \big( |\theta_j(h)| \big).
\end{equation}
Here, $\pi_{\lambda_n, \tau_n} : [0, \infty) \to [0, \infty)$ is a function with two tuning parameters $\lambda_n, \tau_n > 0$, $\theta(h) = \big(\theta_1(h), \cdots, \theta_p(h) \big)^T$ is the parameters of the   the network $h$, and $\theta_j(h)$ is the $j$-th component of $\theta(h)$. 
As in \cite{kurisu2025adaptive}, we assume that $\pi_{\lambda_n, \tau_n}$ satisfies the following conditions:
\begin{itemize}
\item [(i)] $\pi_{\lambda_n, \tau_n}(0) = 0$ and $\pi_{\lambda_n, \tau_n}(\theta)$ is non-decreasing.
\item[(ii)] $\pi_{\lambda_n, \tau_n}(x) = \lambda_n$ if $x > \tau_n$.   
\end{itemize}
As pointed out in \cite{kurisu2025adaptive}, a specific example of the penalty in (\ref{equa_penalty_term_v1}) is the clipped $L_1$ penalty (see \cite{zhang2010analysis}), defined by for all $x \ge 0$ as:
\begin{equation}\label{equa_clipped_L1_penalty}
\pi_{\lambda_n, \tau_n}(x) = \lambda_n \Big(\dfrac{x}{\tau_n} \land 1 \Big), 
\end{equation}
considered also  by \cite{ohn2022nonconvex}, the SCAD penalty considered by \cite{fan2001variable}, the minimax concave penalty see \cite{zhang2010nearly} and the seamless $L_0$ penalty considered in \cite{dicker2013variable}.

\medskip

Numerous contributions have already focused on the NPDNN and/or SPDNN predictors. 
For nonparametric regression from i.i.d. data, \cite{schmidt2020nonparametric} established that the least squares NPDNN estimator can achieve (up to a logarithmic factor), the minimax optimal convergence rate on a class of composition H\"older functions.
In the same vein, \cite{fan2024noise} studied the impact of heavy-tailed noise on the convergence rate of the NPDNN estimator and proved the minimax optimality of this estimator based on the Huber loss, on a class of hierarchical composition functions.
For regression and classification tasks, \cite{ohn2022nonconvex} established that the SPDNN predictor from i.i.d. observations can adaptively attain (up to a logarithmic factor) the minimax optimal rate on classes of smooth and piecewise smooth functions.
\cite{kurisu2025adaptive} focused on nonparametric regression from $\beta$-mixing observations and performed least squares NPDNN and SPDNN estimators.
For the specific case of autoregression, the lower bounds on the class of composition H\"older functions as well as the $\ell^0$-bounded affine class are established. 
In both classes, they proved that the SPDNN estimator adaptively achieves (up to a logarithmic factor) the minimax optimal convergence rate.
For a framework that includes regression and classification, \cite{kengne2025deep} established the convergence rate of the SPDNN predictor, obtained from strongly mixing observations.
In the case of nonparametric autoregression (with the Huber loss) and binary classification (with the logistic loss), they showed that this estimator can adaptively attain (up to a logarithmic factor) the minimax optimal rates on the class of composition H\"older functions.
We also refer to \cite{shen2024nonparametric}, \cite{fan2024factor}, \cite{zhang2024classification},  \cite{kengne2025sparse}, \cite{kengne2025robust} and the references therein, for some other recent advances on NPDNN and SPDNN estimators.
 
 \medskip

The above-mentioned works are performed from i.i.d. data or under some specific dependence conditions.
They do not take into account, for example, $\mathcal{C}$-mixing observations.
Moreover, most of the existing studies are based on a Bernstein-type inequality.
But, as pointed out by \cite{hang2016learning}, they may diverse from one to another since they may be conducted under different assumptions.
One is then interested in considering a general Bernstein-type inequality satisfied by the data and deriving the convergence rate of the DNN estimator.

 \medskip

 In this paper, we develop a unified deep learning setting for various dependent processes, including $\phi$-mixing, $\alpha$-mixing, $\beta$-mixing, and $\mathcal{C}$-mixing processes $\{Z_t = (X_t, Y_t), t\in \Z\}$ which takes values in $\mz = \mx \times \mathcal{Y} \subseteq \R^d \times \R$, based on the training sample $D_n = \{ Z_1 = (X_1, Y_1), \cdots, Z_n =(X_n, Y_n) \}$ and a  loss function $\ell: \R\times \mathcal{Y} \to [0, \infty)$.
We assume that the process $(Z_t)_{t \in \Z}$ fulfills a generalized Bernstein-type inequality (see \cite{hang2016learning} and assumption \textbf{(A4)} below) and our main contributions are as follows:
 
\medskip

\begin{enumerate}
\item \textbf{Oracle inequality for the excess risk of the SPDNN estimator}.
We derive an oracle inequality for the expected excess risk of the SPDNN estimator in a general setting that includes both regression and classification, with a broad class of loss functions.
\item \textbf{Excess risk bounds.}
Excess risk bounds of the proposed DNN estimators are derived on the classes of H\"older and composition H\"older smooth functions.
\begin{itemize}
\item When the target function belongs to the class of H\"older functions with a smoothness parameter $s>0$, it is established that the convergence rate of the expected excess risk of the NPDNN and SPDNN estimators is of order $\mathcal{O}\Bigg( \Big(\varphi(n)\Big)^{-\frac{\kappa s}{\kappa s + d}} \Big(\log \big(\varphi(n)\big)\Big)^ 3 \Bigg)$ (where $\varphi(n)$ is given in assumption \textbf{(A4)} for some constant $\kappa \ge 1$, and $d$ is the dimension of the input space).
\item On the class of composition H\"older functions (defined in Subsection \ref{excess_risk_comp_Holder}, see (\ref{equa_compo_structured_function})),
the convergence rate above is of order $\mathcal{O} \Bigg( \Big(\phi _{n, \varphi} ^{\kappa/2} \lor \phi _{n, \varphi}\Big) \Big(\log \big( \varphi(n) \big)  \Big)^{3} \Bigg)$, where $ \phi _{n, \varphi}$ is defined in (\ref{def_phi_nalpha}). 
\end{itemize}
\item \textbf{Applications and examples.} 
Applications to i.i.d. observations, $\phi$-mixing, $\alpha$-mixing, and $\mathcal{C}$-mixing processes are carried out.
For each of these examples, the generalized Bernstein-type inequality in Assumption \textbf{(A4)} holds and the bounds above are applied with $\varphi(n)$ as in Table \ref{Appl_example}.

\begin{table}[h!]
%\footnotesize
%\scriptsize
%\footnotesize
\centering
\caption{ \it Applications and examples of some dependence structures. The convergence rate is for the NPDNN and SPDNN estimators. The classes $\mathcal{C}^s(\mathcal{X}, \mathcal{K})$ and $\mathcal{G}(q, \bold{d}, \bold{t}, \boldsymbol{\beta}, A)$ are defined in (\ref{equa_ball_Holder}) and (\ref{equa_compo_structured_function}). $\kappa$ is given in Assumption \textbf{(A3)} and $\varrho$ is given in (\ref{expo_alpha_mixing_process}), (\ref{equa_phi_mixing_coef}) and (\ref{equa_poly_C-mixing_coef}) .}
\label{Appl_example}
\vspace{.1cm}
\begin{tabular}{cccc}
\hline  
  &  &  &    \\  
  Dependence structure & $\varphi(n)$ & Convergence rate on $\mathcal{C}^s(\mathcal{X}, \mathcal{K})$ & Convergence rate on $\mathcal{G}(q, \bold{d}, \bold{t}, \boldsymbol{\beta}, A)$   \\
   &  &  &   \\  
\Xhline{.6pt}%\\                  
\rule[0cm]{0cm}{.35cm}
 I.I.d. observations  & $n$  & $n^{-\frac{\kappa s}{\kappa s + d}} \big(\log n\big)^3$ & $\big(\phi _{n} ^{\kappa/2} \lor \phi _{n}\big) \big(\log n \big)^{3}$   \\
           &  &  &    \\       
  $\phi$-mixing  & $n$  & $n^{-\frac{\kappa s}{\kappa s + d}} \big(\log n\big)^3$ & $\big(\phi _{n} ^{\kappa/2} \lor \phi _{n}\big) \big(\log n \big)^{3}$   \\
             &  &  &    \\       
 Exponential $\alpha$-mixing  & $n / (\log n)^2$ & $n^{-\frac{\kappa s}{\kappa s + d}} \big(\log n\big)^5$  & $\big(\phi _{n} ^{\kappa/2} \lor \phi _{n}\big) \big(\log n \big)^{3 + 2 \lor \kappa}$\\ 
            &  &  &    \\       
   Subexponential $\alpha$-mixing  & $n^{\varrho /(\varrho + 1)}$ & $n^{- \frac{\varrho}{\varrho + 1} \frac{\kappa s}{\kappa s + d}} \big(\log n\big)^3$  & $  \big(\phi _{n} ^{\frac{\kappa \varrho}{2(\varrho + 1)}} \lor \phi _{n}^{\frac{\varrho}{\varrho + 1}} \big) (\log n)^3  $  \\ 
            &  &  &    \\       
    Geometrically $\mathcal{C}$-mixing &  $n/\big(\log n\big)^ {2/\varrho} $  & $n^{-\frac{\kappa s}{\kappa s + d}} \big(\log n \big)^ {3 + 2/\varrho}$ & $ \big(\phi_{n}^{\kappa/2} \lor \phi_{n} \big)(\log n)^ {3 + (2/\varrho) \lor \kappa}$  \\
            &  &  &   \\       
   Polynomially $\mathcal{C}$-mixing  &  $ n^{(\varrho - 2)/(\varrho + 1)} $  & $ n^{-\frac{\varrho - 2}{\varrho + 1}\frac{\kappa s}{\kappa s + d}} \big(\log n \big)^3 $&  $  \big(\phi _{n} ^{\frac{\kappa (\varrho - 2)}{2(\varrho + 1)}} \lor \phi _{n}^{\frac{\varrho-2}{\varrho + 1}} \big) (\log n)^3  $  \\ 
 \Xhline{.9pt}
 %\hline
\end{tabular}
%\scriptsize
\end{table} 

For these dependence structures, the NPDNN and SPDNN estimators can achieve (up to a logarithmic factor), the minimax optimal convergence rate.
For example, in the case of the regression with a symmetric (around 0) error and the Huber loss, and the classification with the logistic loss, Assumption \textbf{(A3)} holds with $\kappa = 2$.
In these cases and for i.i.d. data, $\phi$-mixing, exponential $\alpha$-mixing and geometrically $\mathcal{C}$-mixing processes, the convergence rates (up to a logarithmic factor) on the classes of H\"older functions $\mathcal{C}^s(\mathcal{X}, \mathcal{K})$ and composition H\"older functions $\mathcal{G}(q, \bold{d}, \bold{t}, \boldsymbol{\beta}, A)$ are $n^{-\frac{2 s}{2 s + d}} $ and $\phi_n$ respectively. These rates are optimal (up to a logarithmic factor) in the minimax sense.  
\end{enumerate} 
\noindent
In comparison with \cite{kengne2025robust}, the results in Subsection \ref{excess_risk_bound_NPDNN} below coincide with well known results for classical models and can reach (up to a logarithmic factor)  the minimax optimal rate, which is not the case for the bounds in \cite{kengne2025robust}. 
The theoretical results on the SPDNN estimator in Subsection \ref{susec_spdnn} are extensions of \cite{kengne2025deep} to a broad class of processes, including $\mathcal{C}$-mixing observations. 
 
\medskip

The rest of the paper is outlined as follows.
Section \ref{asump} presents some notations and assumptions. 
In Section \ref{oracle_inqua_excess_risk}, we establish an oracle inequality for the SPDNN predictor and derive excess risk bounds of the NPDNN and  SPDNN estimators.
Section \ref{application_examples} focuses on applications and examples for some processes that satisfy the generalized Bernstein-type inequality.
Example of nonparametric autoregression with exogenous covariate is considered in Section \ref{sect_np_reg}.
Section \ref{prove} is devoted to the proofs of the main results.

\medskip

\section{Notations and assumptions}\label{asump}
\subsection{Some notations}
Let $E_1$ and $E_2$ be two subsets of separable Banach spaces equipped respectively with a norm $\|\cdot\|_{E_1}, \|\cdot\|_{E_2}$.
In the sequel, the following notations will be used.  
\begin{itemize}
\item For $h: E_1 \rightarrow E_2$ and $U \subseteq E_1$, we set,
\begin{equation}\label{def_norm_inf}
\| h\|_{\infty} = \sup_{x \in E_1} \| h(x) \|_{E_2}, ~ \| h\|_{\infty,U} = \sup_{x \in U} \| h(x) \|_{E_2}.
\end{equation}   
\item For any $h: E_1 \rightarrow E_2, \epsilon >0$, set, 
\[ B (h, \epsilon) = \big\{f: E_1 \rightarrow E_2, ~ \| f - h\|_\infty \leq \epsilon \big\},  \]
where $\| \cdot \|_\infty$ denotes the sup-norm defined in (\ref{def_norm_inf}).
\item Let $\mathcal{H} \subset \mathcal{F}(E_1, E_2)$ (the set of measurable functions from $E_1$ to $E_2$).
 For any $\epsilon > 0$, the $\epsilon$-covering number $\mathcal{N}(\mathcal{H}, \epsilon)$  of $\mathcal{H}$, is given by,
\begin{equation}\label{epsi_covering_number}
 \mathcal{N}(\mathcal{H}, \epsilon) = \inf\Big\{ m \geq 1 ~: \exists h_1, \cdots, h_m \in \mathcal{H} ~ \text{such that} ~ \mathcal{H} \subset \bigcup_{i=1}^m B(h_i, \epsilon) \Big\}.
\end{equation}
It is the minimal number of balls of radius $\epsilon$ needed to cover $\mathcal{H}$. 
\item For any $x \in \R$, $\lfloor x \rfloor$ denotes the largest integer strictly smaller than $x$.
\item For all $u, v \in \R$, we set $u\lor v = \max(u, v)$ and $u \land v = \min(u, v)$.
\item For two sequences of real numbers $(u_n)$ and $(v_n)$, we write $ u_n \lesssim v_n$ or $ v_n \gtrsim u_n$ if there exists a constant $C > 0$ such that $u_n \leq  C v_n$ for all $n \in \N$; $u_n \asymp v_n$ if $u_n \lesssim v_n$ and $u_n \gtrsim v_n$.
\end{itemize}

\medskip

\noindent
Let us consider the following definition of  piecewise linear and  locally quadratic functions, see also \cite{ohn2019smooth}, \cite{ohn2022nonconvex}.
\begin{Def}\label{def_pwl_quad}
Let a function $g: \R \rightarrow \R$.
\begin{enumerate}
\item $g$ is continuous piecewise linear (or ``piecewise linear" for notational simplicity) if it is continuous and there exists $K$ ($K\in \N$) break points $a_1,\cdots, a_K \in \R$ with $a_1 \leq a_2\leq\cdots \leq a_K $ such that, for any $k=1,\cdots,K$, $g'(a_k-) \neq g'(a_k+)$ and $g$ is linear on $(-\infty,a_1], [a_1,a_2],\cdots [a_K,\infty)$.
\item $g$ is locally quadratic if there exists an interval $(a,b)$ on which $g$ is three times continuously differentiable with bounded derivatives and there exists $t \in (a,b)$ such that $g'(t) \neq 0$ and $g''(t) \neq 0$.
\end{enumerate} 
\end{Def}
\subsection{Some assumptions}\label{subsec_assump}
We consider a process $\{ Z_t = (X_t, Y_t), t \in \Z \}$ which take values in $ \mathcal{Z}= \mx \times \mathcal{Y} \subset \R^d \times \R$, a loss function 
$\ell : \R \times \mathcal{Y} \to [0, \infty)$, an activation function $\sigma: \R \rightarrow \R $, and  set the following assumptions.
\begin{itemize}
\item[\textbf{(A0)}:] $\mathcal{X}\subset R^d$ is a compact set.
\item[\textbf{(A1)}:] There exists a constant $C_{\sigma}>0$ such that the activation function $\sigma$ is $C_{\sigma}$-Lipschitz, that is, there exists $C_{\sigma} > 0$ such that $|\sigma(x_1) - \sigma(x_2)| \leq C_{\sigma} |x_1 - x_2|$ for any $x_1, x_2 \in \R$.
Moreover, $\sigma$ is  either piecewise linear or locally quadratic and fixes a non empty interior segment $I \subseteq [0,1]$ (i.e. $\sigma(z) = z$ for all $z \in I$).
\item[\textbf{(A2)}:] There exists $\mk _{\ell} > 0$ such that, for any $(y_1, y_2) \in \R ^2$, we have $|\ell(y_1, y) - \ell(y_2, y)| \leq \mk _{\ell} |y_1 - y_2|$ for all $y \in \mathcal{Y}$.
\item[\textbf{(A3)}:] 
Local structure of the excess risk: There exist three constants $\mk_0:= \mk_0(Z_0, \ell, h^*) , \varepsilon_0:= \varepsilon_0(Z_0, \ell, h^*) >0$ and $\kappa:=\kappa(Z_0, \ell, h^*) \geq 1$ such that,  
\begin{equation}\label{assump_local_quadr}
 R(h) - R(h^*) \leq \mk_0 \| h - h^*\|^\kappa_{\kappa, P_{X_0}}, 
\end{equation}
for any measurable function $h: \R^d \rightarrow \R$ satisfying $\|h - h^*\|_{\infty, \mx} \leq \varepsilon_0$; where $P_{X_0}$ denotes the distribution of $X_0$ and 
\[ \| h - h^{*}\|_{r, P_{X_0}} ^r := \displaystyle \int \| h (\text{x})- h^{*} (\text{x}) \|^r  d P_{X_0} ( \text{x}),  
\]
for all $r \geq 1$.
\item [\textbf{(A4)}:] The $\mz$-valued process $\{ Z_t = (X_t, Y_t), t \in \Z \}$ is stationary and ergodic.
 Moreover, let $g: \mz \to \R$ be a bounded measurable function such that $\E[g(Z_1)] = 0$, $ \E[\big( g(Z_1) \big) ^2] \leq \gamma ^2$ and $\|g\|_{\infty} \leq A$ for some $A > 0$ and $\gamma \ge 0$. Assume that, for all $\varepsilon > 0$, there exist $n_0 \ge 1$ independent of $\varepsilon$ and $\varphi(n) \ge 1$ such that for all $n \ge n_0$, we get
\begin{equation}\label{G_beernstein_expo_inequa}
P\Bigg\{\dfrac{1}{n} \sum _{i=1} ^n g(Z_i) \ge \varepsilon \Bigg\} \leq C\exp \Bigg( -\dfrac{\varepsilon ^2 \varphi(n)}{c_{\gamma} \gamma ^2 + c_A \varepsilon A}\Bigg),    
\end{equation}
where $\varphi: \N \to \N$ is a non-decreasing function, satisfying $\varphi(n) \leq n$, $C>0$ is a constant independent of $n$, and $c _{\gamma}, c _A$ are positive constants.

\end{itemize}

\medskip

\noindent
The assumptions above are satisfied in several classical frameworks.
 The ReLU (rectified linear unit) activation function defined by $\sigma(x) = \max(x, 0)$, satisfies Assumption \textbf{(A1)}. 
We also refer to \cite{kengne2025excess} for several other activation functions that fulfill \textbf{(A1)}. 
  For the $L_1$ and the Huber (with parameter $\delta >0$) losses, \textbf{(A2)} is satisfied with $\mathcal{K}_{\ell} = 1$ and $\mathcal{K}_{\ell} = \delta$ respectively. 
For all $y \in \R$ and $y' \in \mathcal{Y}$, recall:
\begin{itemize}
\item The $L_1$ loss: $\ell(y,y') = |y-y'|$;
\item The Huber loss with parameter $\delta >0$:
   \[   \ell_\delta(y,y') =  \left\{
\begin{array}{ll}
      \frac{1}{2}(y-y')^2 &  \text{ if } ~ |y-y'|\leq \delta , \\
 \delta |y-y'| - \frac{1}{2} \delta^2 &   ~ \text{ otherwise}.
\end{array}
\right.    \]
\end{itemize}
 Assumption \textbf{(A3)} is a local condition on the excess risk. When the loss function is Lipschitz continuous, this condition is satisfied with $\kappa =1$. 
For nonparametric regression with a symmetric (around 0) error and with the Huber loss, one can get from Proposition 3.1 in \cite{fan2024noise} that \textbf{(A3)} is satisfied with $\kappa=2$.
In the case of binary classification with the logistic loss, one can get from Lemma C.6 in \cite{zhang2024classification} (see also \cite{alquier2025minimax}) that Assumption \textbf{(A3)} holds with $\kappa = 2$.
\textbf{(A4)} is satisfied for several classical processes, see Section \ref{application_examples}.

\section{Excess risk bounds}\label{oracle_inqua_excess_risk}

\subsection{Excess risk bounds of the NPDNN estimator}\label{excess_risk_bound_NPDNN}
In this subsection, we establish upper bounds of the excess risk of the NPDNN estimator for observations that satisfy the generalized Bernstein-type inequality in Assumption \textbf{(A4)}. 
In the next theorem, the convergence rate of the excess risk is derived on the class of  H\"older smooth functions.
Let us recall:

\medskip
\noindent
\textbf{Class of H\"older smooth functions.}\label{excess_risk_Holder_space}
 Recall that for any $r\in \N, D\subset \R^r, \beta, A >0$, the ball of $\beta$-H\"older functions with radius $A$ is defined by:
\begin{equation}\label{equa_ball_Holder}
 \mathcal{C}^{\beta}(D, A)= \Bigg\{h : D \rightarrow \R : \underset{\boldsymbol{\alpha} \in \N^r : |\boldsymbol{\alpha}|_1 < \beta}{\sum}\|\partial^{\boldsymbol{\alpha}}h\|_{\infty} + \underset{\boldsymbol{\alpha} \in \N^r : |\boldsymbol{\alpha}|_1 = \lfloor \beta \rfloor}{\sum}~\underset{x\ne y}{\underset{x, y\in D}{\sup}} \frac{|\partial^{\boldsymbol{\alpha}}h(x) - \partial^{\boldsymbol{\alpha}}h(y)|}{|x - y|^{\beta - \lfloor \beta \rfloor}} \leq A \Bigg\},  
\end{equation}
with $\boldsymbol{\alpha} = (\alpha_1, \dots, \alpha_r)\in \N^r$, $|\boldsymbol{\alpha}|_1 \coloneqq \sum_{i=1}^r \alpha_i$ and $\partial^{\boldsymbol{\alpha}} = \partial^{\alpha_1} \dots \partial^{\alpha_r}$.

\begin{thm} \label{thm2} 
Assume (\textbf{A0}), (\textbf{A1}), (\textbf{A2}), \textbf{(A4)}. Assume that \textbf{(A3)} is satisfied for some $\kappa \geq 1$, $\mk_0, \varepsilon_0>0$ and that  $h^{*} \in \mathcal{C}^{s}(\mathcal{X}, \mathcal{K})$ for some $s, \mathcal{K} > 0$, where $h^{*}$ is defined in (\ref{best_pred_F}). Set 
$L_n = \dfrac{s L_0}{\kappa s + d} \log\big(\varphi(n)\big), N_n = N_0 \big(\varphi(n)\big)^{\frac{d}{ \kappa s + d}}, S_n = \frac{s S_0}{\kappa s + d} \big(\varphi(n)\big)^{\frac{d}{\kappa s + d}}\log\big(\varphi(n)\big), 
$
 $ B_n = B_0 \big(\varphi(n)\big)^{\frac{4(d + s)}{\kappa s + d}} $ and $F_n = F$, 
for some positive constants $L_0, N_0, S_0, B_0 > 0$. Consider the class of DNN $\mathcal{H}_{\sigma}(L_n, N_n, B_n, F_n, S_n)$ defined in (\ref{constrained_dnn_class}).
Then, there exists $n_0 = n_0(\mk_{\ell}, \mk, F, \kappa, L_0, N_0, B_0, S_0, C_\sigma,$ $ s, d,) \in \N$, such that for all $n \geq n_0$, the NPDNN estimator $\widehat{h}_{n, NP}$ defined in  (\ref{NP_DNNs_Estimators}) satisfies,
\begin{equation}\label{thm2_excess_risk_bound}
\E[R(\widehat{h}_{n, NP}) - R(h^{*})] \lesssim \dfrac{ \big(\log \big(\varphi(n)\big)  \big)^{\nu}}{ \big(\varphi(n)\big)^{\frac{\kappa s}{\kappa s + d}} },
\end{equation}
 for all $\nu>3$.
\end{thm}

\medskip

In the sequel, we also deal with the class of composition H\"older functions.

\noindent
\textbf{Class of composition H\"older functions.}\label{excess_risk_comp_Holder}
Let $q\in \N, \boldsymbol{d} = (d_0, \dots, d_{q+1})\in \N^{q + 2}$ with $d_0 = d, d_{q+1} = 1, \boldsymbol{t} = (t_0, \dots, t_q)\in \N^{q+1}$ with $t_i \leq d_i$ for all $i$, $\boldsymbol{\beta} = (\beta_0, \dots, \beta_q)\in (0,\infty)^{q+1}$ and $A>0$. 
For all $l < u$, denote by $\mathcal{C}_{t_i}^{\beta_i}([l, u]^{d_i}, A)$ the set of functions $f: [l, u]^{d_i} \rightarrow \R$ that depend on at most $t_i$ coordinates and $f \in \mathcal{C}^{\beta_i}([l, u]^{d_i}, A)$.

\medskip

\noindent Define the set of composition H\"older smooth functions $\mathcal{G}(q, \bold{d}, \bold{t}, \boldsymbol{\beta}, A)$ studied in \cite{schmidt2020nonparametric}, by: 
\begin{multline}\label{equa_compo_structured_function}
\mathcal{G}(q, \bold{d}, \bold{t}, \boldsymbol{\beta}, A) \coloneqq \Big\{h = g_q\circ \dots \circ g_0, g_i = (g_{ij})_{j=1,\cdots,d_{i+1}} : [l_i, u_i]^{d_i}
 \rightarrow [l_{i+1}, u_{i+1}]^{d_{i+1}}, g_{ij}\in C_{t_i}^{\beta_i}
 ([l_i, u_i]^{d_i}, A) ~ \\
        \text{for some } l_i, u_i \in \R \text{ such that } ~ |l_i|, |u_i| \leq A, \text{ for } i=1,\cdots,q \Big\}.
\end{multline}
The smoothness of a function in the class $\mathcal{G}(q, \bold{d}, \bold{t}, \boldsymbol{\beta}, A)$ (see (\cite{juditsky2009nonparametric})) is given by
$\beta_i^* \coloneqq \beta_i\prod_{j = i+1}^q(\beta_{j} \land 1)$.
 Under Assumption \textbf{(A4)}, we set:
\begin{equation}\label{def_phi_nalpha}
\phi_{n, \varphi} \coloneqq  \underset{0\leq i \leq q}{\max}\big(\varphi(n) \big)^{ \frac{-2\beta_i^* }{2\beta_i^* + t_i}   }.
\end{equation}

\noindent
\medskip

\noindent
The next Theorem is established with the ReLU activation function, that is $\sigma(x) = \max(x, 0)$ (that satisfies \textbf{(A1)} with $C_\sigma =1$) and $\mathcal{X} = [0,1]^d$.
\begin{thm}\label{thm3}
Set $\mathcal{X} = [0,1]^d$. Assume that \textbf{(A2)}, \textbf{(A4)} and  \textbf{(A3)} for some $\kappa \ge 1, \mk_0 , \varepsilon_0 > 0$ hold and that $h ^* \in \mathcal{G}(q, \bold{d}, \bold{t}, \boldsymbol{\beta}, \mathcal{K})$. Consider the DNN class $\mathcal{H}_{\sigma}(L_n, N_n, B_n, F_n, S_n)$ with network architecture $(L_n, N_n, B_n, F_n, S_n)$ and an activation function $\sigma$. 
Let $L_n\asymp \log\big( \varphi(n) \big), N_n \asymp \varphi(n) \phi_{n, \varphi}, B_n = B \geq 1 , F_n = F > \max(\mathcal{K}, 1), S_n \asymp  \varphi(n) \phi_{n, \varphi} \log \big( \varphi(n) \big)$.
Then, the NPDNN estimator defined in (\ref{NP_DNNs_Estimators}) with the ReLU activation function satisfies,  
\begin{equation}\label{them_comp_npdnn}
\E[R(\widehat{h}_{n,NP}) - R(h^{*})] \leq  \Big(\phi _{n, \varphi} ^{\kappa/2} \lor \phi _{n, \varphi}\Big) \Big(\log \big( \varphi(n) \big)  \Big)^{\nu}, 
\end{equation}
for all $\nu > 3$.
\end{thm}
\noindent
When $\kappa=2$, the bounds in (\ref{thm2_excess_risk_bound}) and (\ref{them_comp_npdnn}) are optimal (up to a logarithmic factor) for several classical models, see Section \ref{application_examples} below.

\subsection{Excess risk bounds of the SPDNN estimator}\label{susec_spdnn}
This subsection deals with the SPDNN estimator.
The results of this subsection are extensions of the results of \cite{kengne2025deep} to processes that satisfy the generalized Bernstein-type inequality in Assumption \textbf{(A4)}, including $\mathcal{C}$-mixing processes. 
In the sequel, we also need to calibrate the tuning parameters $\lambda_n, \tau_n$ of the penalty term, in addition to the network architecture parameters $L_n, N_n, B_n, F_n$. 
The following theorem establishes an oracle inequality of the expected excess risk for the SPDNN estimator.
\begin{thm}\label{thm1}
Assume that \textbf{(A1)}, \textbf{(A2)}, \textbf{(A4)}  hold.  Let $ L_n \asymp \log \big( \varphi(n) \big), N_n \lesssim \big( \varphi(n) \big)^{\nu_1}, 1 \leq B_n \lesssim \big( \varphi(n) \big)^{\nu_2}, F_n = F > 0$, for some $ \nu_1 > 0, \nu_2 > 0$; and that the target function $h^{*}$ defined in (\ref{best_pred_F}) is such that $\|h^{*}\|_{\infty} \leq \mathcal{K}$ for some constant $\mathcal{K} > 0$. Then, there exists $n_0 \in \N$ such that, for all $n \geq n_0$, the SPDNN estimator $\widehat{h}_{n,SP}$ defined in (\ref{sparse_DNNs_Estimators}), 
with $\lambda_n \asymp \frac{\big( \log(\varphi(n)) \big) ^{\nu_3}}{\varphi(n)}$ for some $\nu_3>2$ and $\tau_n \leq \dfrac{1}{16 \mathcal{K}_{\ell}(L_n + 1)((N_n + 1)B_n)^{L_n + 1} \varphi(n) }$ satisfies,
\begin{align}
\nonumber \E[R(\widehat{h}_{n,SP}) - R(h^{*})]
  \leq 2\Big(\underset{h \in \mathcal{H}_{\sigma}(L_n, N_n, B_n, F)} {\inf} R(h) - R(h^{*}) + J_n(h) \Big]  \Big)  + \dfrac{\Xi}{\varphi(n)},
\end{align}
for some constant $\Xi>0$, and $ \mathcal{K}_{\ell} > 0$ is given in Assumption (\textbf{A2}).
\end{thm}

\medskip

\noindent
The following corollary derives upper bounds of the expected excess risk of the SPDNN estimator on the classes of H\"older and composition H\"older functions.
This corollary can be obtained in the same way as in the proofs of Corollary 4.3. and Corollary 4.5 in \cite{kengne2025deep}. The proof is omitted.

\begin{Corol}\label{corol1}
Assume \textbf{(A2)}, \textbf{(A3)} for some $\kappa \ge 1, \mk_0, \varepsilon_0 > 0$ and \textbf{(A4)}. 
Consider the DNN class $\mathcal{H}_{\sigma}(L_n, N_n, B_n, F_n)$ with network architecture $(L_n, N_n, B_n, F_n)$ and the SPDNN estimator $\widehat{h}_{n,SP}$ defined in (\ref{sparse_DNNs_Estimators}) where  $\lambda_n, \tau_n$ are given as in Theorem \ref{thm1}.
\begin{enumerate}
\item If \textbf{(A0)}, \textbf{(A1)} hold and $L_n \asymp \log\big( \varphi(n) \big), N_n \asymp \big(\varphi(n) \big)^{\frac{d}{\kappa s + d}}, B_n \asymp \big(\varphi(n) \big)^{\frac{4(s + d)}{\kappa s + d}}, F_n = F >0$, then 
  \begin{equation}\label{excess_risk_bound_v1}
\underset{h^{*} \in  \mathcal{C}^{s}(\mathcal{X}, \mathcal{K})}{\sup} \big( \E[R(\widehat{h}_{n,SP}) - R(h^{*})] \big) \lesssim \dfrac{  \big(\log\big( \varphi(n) \big) \big)^{\nu} }{\big(\varphi(n)\big)^{\kappa s/(\kappa s + d)}},
\end{equation}
for all $s, \mk >0$ and $\nu > 3$.
\item Assume that $\mathcal{X} = [0,1]^d$ and that $\sigma$ is the ReLU activation function and consider the class of composition H\"older functions $\mathcal{G}(q, \bold{d}, \bold{t}, \boldsymbol{\beta}, \mathcal{K})$ defined in (\ref{equa_compo_structured_function}).
If $L_n\asymp \log\big( \varphi(n) \big), N_n \asymp \varphi(n) \phi_{n, \varphi}, B_n = B \geq 1 , F_n = F > \max(\mathcal{K}, 1)$, then
\begin{align}\label{excess_risk_bound_corol01_v1}
\underset{h^{*} \in \mathcal{G}(q, \bold{d}, \bold{t}, \boldsymbol{\beta}, \mathcal{K})}{\sup} \big( \E[R(\widehat{h}_{n,SP}) - R(h^{*})] \big)  & \lesssim \Big( \phi_{n, \varphi} ^{\kappa/2} \lor \phi_{n, \varphi} \Big) \big(\log\big( \varphi(n) \big) \big)^{\nu},  
\end{align}
for all $\nu > 3$.
\end{enumerate}
\end{Corol}
The next section performs applications and examples to some well known dependence structures. 

\section{Applications and examples}\label{application_examples}
Let us consider some examples to which the theoretical results above can be applied.
Recall that (see Subsection \ref{subsec_assump}) for a regression problem with the Huber loss and symmetric errors around $0$, as well as for a classification task with logistic loss, assumption \textbf{(A3)} holds with $\kappa = 2$. 
In the sequel, we consider the class of composition H\"older functions $\mathcal{G}(q, \bold{d}, \bold{t}, \boldsymbol{\beta}, A)$ defined in (\ref{equa_compo_structured_function}) and set 
\begin{equation}\label{def_phi_n_alpha}
\phi_n \coloneqq  \underset{0\leq i \leq q}{\max} n^{ \frac{-2\beta_i^* }{2\beta_i^* + t_i}   },
\end{equation}
$\beta_i^* \coloneqq \beta_i\prod_{j = i+1}^q(\beta_{j} \land 1)$ is the smoothness of a composition function in $\mathcal{G}(q, \bold{d}, \bold{t}, \boldsymbol{\beta}, A)$ (see for instance (\cite{juditsky2009nonparametric})).
We will also deal with the class of H\"older smooth functions $\mathcal{C}^{s}(\mathcal{X}, \mathcal{K})$ for some $s, \mathcal{K} > 0$.
Furthermore, for any stationary process $Z = \{Z_t\}_{t\in \Z}$ on a probability space $(\Omega, \mathcal{B}, P)$, we denote by $\sigma_{t}^{\infty}$ and $\sigma_{-\infty}^t$ for  $-\infty \leq t \leq \infty$, the $\sigma$-algebras generated respectively by the random variables $Z_j, j \ge t$ and $Z_j, j\leq t$.

\subsection{I.I.d. processes.}
Assume that the process $\{(X_t, Y_t), t \in \Z \} $ is i.i.d..
The inequality (\ref{G_beernstein_expo_inequa}) in assumption \textbf{(A4)} holds with $n_0 = 1$, $C = 1$, $c_\gamma = 2$, $c_A = \dfrac{2}{3}$, and $\varphi(n) = n$.
Thus, the theoretical results above can be applied. Therefore:
\begin{itemize}
\item The convergence rate of the expected excess risk of the NPDNN estimator on the class of H\"older and composition H\"older functions are $\mathcal{O} \Big(n^{-\frac{\kappa s}{\kappa s + d}}  (\log n)^ 3 \Big)$ and $\mathcal{O} \Big( \big(\phi_n^{\kappa/2} \lor \phi_n \big)(\log n)^ 3\Big)$ respectively, with $\kappa \ge 1$ given in \textbf{(A3)}.
\item For the SPDNN predictor, the same rates as above are obtained.
\end{itemize}

\medskip

\noindent
When $\kappa = 2$, the convergence rate on the class of composition H\"older functions, is of order $\mathcal{O} \Big(\phi_n (\log n)^ 3\Big)$ for both NPDNN and SPDNN estimators. 
This rate matches (up to a logarithmic factor) with the lower bound established by \cite{schmidt2020nonparametric} and \cite{alquier2025minimax} for regression and classification.
In this case, the NPDNN and the SPDNN estimators proposed are optimal (up to a logarithmic factor) in the minimax sense.

\subsection{$\phi$-mixing processes}. 
A $\mathcal{Z}$-valued process $Z = \{Z_t\}_{t\in \Z} $ is said to be $\phi$-mixing, if it satisfies
\[\underset{A \in \sigma_{-\infty}^0, B \in\sigma_k^{\infty}}{\sup} \{ |P(B) - P(B|A|\} = \phi( k ) \to 0 \text{ as } k \to \infty.
\]
where $\phi(k)$ is called the $\phi$-mixing coefficient.

\medskip

\noindent
Under $\phi$-mixing condition on the process $\{(X_t, Y_t), t \in \Z \} $, the inequality (\ref{G_beernstein_expo_inequa}) in assumption \textbf{(A4)} fulfills (see \cite{samson2000concentration}) with $n_0 = 1$, $C = 1$, $c_\gamma = 32 \sum_{k \geq 1} \sqrt{\phi(k)}$, $c_A = 8\sum_{k \geq 1} \sqrt{\phi(k)}$, and $\varphi(n) = n$. 
So, the theoretical results established can be applied.
 Also, the same rates and conclusion as in the i.i.d. case are obtained.

\subsection{Strong mixing processes}
A $\mathcal{Z}$-valued process $Z = \{Z_t\}_{t\in \Z} $ is said to be $\alpha$-mixing or strongly mixing if it satisfies
\[\underset{A \in \sigma_{-\infty}^0, B \in\sigma_k^{\infty}}{\sup} \{ |P(A \cap B) - P(A)P(B) |\} = \alpha(k ) \to 0 \text{ as } k \to \infty.
\]
where $\alpha( k)$ is called the $\alpha$-mixing coefficient.
\begin{enumerate}
\item  \textbf{Exponential $\alpha$-mixing processes.} 
The stochastic process $Z$ is said to be exponential strongly mixing, if the mixing coefficient $\alpha(k)$ satisfies
\begin{equation}\label{expo_alpha_mixing_process}
\alpha(k) \leq c\exp \big(-bk^{\varrho}\big),  ~~ \text{ for all } k \ge 1,  \end{equation}
for some constants $b > 0$, $c \geq 0$, and $\varrho \ge 1$.
\medskip

\noindent
From \cite{merlevede2009bernstein} and \cite{hang2016learning}, one can find $n_0 \geq 1$, $C, c_\gamma, c_A \geq 0$ such that, \textbf{(A4)} holds for all $n \geq n_0$ with $\varphi(n) = n/\big(\log n \big)^2$.
Thus, the theoretical studies above applied to exponential strong mixing processes give for $\kappa \geq 1$ given in \textbf{(A3)}:
\begin{itemize}
\item A convergence rate of order $\mathcal{O} \Big(n^{-\frac{\kappa s}{\kappa s + d}} \big( \log (n) \big)^5 \Big)$ on the class of H\"older smooth functions, for the NPDNN and the SPDNN predictors.
\item This rate is  $\mathcal{O} \Big( \big(\phi_n^{\kappa/2} \lor \phi_n \big)(\log n)^{3 + 2 \lor \kappa}\Big)$ on the class of composition H\"older functions. 
\end{itemize}
When $\kappa=2$, the same conclusion as in the i.i.d. case regarding the minimax optimality is applied.
\item \textbf{Subexponential $\alpha$-mixing processes.}
The stochastic process $Z$ is said to be subexponential strongly mixing, if the mixing coefficient $\alpha(k)$ satisfies (\ref{expo_alpha_mixing_process}) for some $\varrho > 0$.

\medskip

\noindent
From \cite{modha1996minimum} and \cite{hang2016learning}, we get that the inequality (\ref{G_beernstein_expo_inequa}) in assumption \textbf{(A4)} holds with $n_0 = \max\big( b/8, 2^{2 + 5/\varrho}b^{-1/\varrho}\big)$, $C = 1 + 4e^{-2}c$, $c_ {\gamma} = \big(8^{2 + \varrho/b}\big)^ {\frac{1}{1 + \varrho}}$, $c_ A = \dfrac{\big(8^{2 + \varrho/b}\big)^ {\frac{1}{1 + \varrho}} }{3}$, and $\varphi(n) = n^{\frac{\varrho}{\varrho + 1}}$. 
Thus, the theoretical results above can be applied. Hence, we obtain for $\kappa \geq 1$ given in \textbf{(A3)}:
\begin{itemize}
\item When the target function belongs to the class of H\"older smooth functions, the convergence rate of the NPDNN and the SPDNN estimators is of order $\mathcal{O} \Big(n^{-\frac{\varrho}{\varrho + 1}\frac{\kappa s}{\kappa s + d}} \big(\log n \big)^ 3 \Big)$.
\item This rate is of order  $\mathcal{O} \Big( \big(\phi _{n} ^{\frac{\kappa \varrho}{2(\varrho + 1)}} \lor \phi _{n}^{\frac{\varrho}{\varrho + 1}} \big) (\log n)^3 \Big)$ with $\kappa \ge 1$, when we deal with the class of composition H\"older functions. 
\end{itemize}
\end{enumerate}

\subsection{$\mathcal{C}$-mixing processes}
 Recall that $(\Omega, \mathcal{B}, P)$ denotes a probability space. For any $p \geq 1$, let $L_p(P) := \{f : \Omega \to \R \text{ measurable: } \|f\|_p < \infty\}$, where $\|f\|_p = \big( \int_\Omega |f(x)| dP(x)  \big)^{1/p} $. In addition, if $\mathcal{B}' \subset \mathcal{B}$ is a sub-$\sigma$-algebra, we denote by $L_p(\mathcal{B'}, P) := \{f : \Omega \to \R,  ~\mathcal{B'}-\text{measurable: } \|f\|_{p} < \infty \}$.
 Furthermore, given a semi-norm $\| \cdot\|$ on the vector space of bounded measurable functions $f : \mathcal{Z} \rightarrow \R$ ($\mathcal{Z} = \mathcal{X}  \times \mathcal{Y} \subset \R^d \times \R $), define the $\mathcal{C}$-norm by
$\|f\|_{\mathcal{C}} := \|f\|_{\infty} + \|f\|$, (where $\|\cdot\|_{\infty}$ is the sup-norm defined in (\ref{def_norm_inf})) and denote the set of all bounded $\mathcal{C}$-functions by $\mathcal{C}(\mathcal{Z}) := \{f : \mathcal{Z} \to \R: \|f\|_{\mathcal{C}} < \infty\}$.  

\medskip
\noindent
A $\mz$-valued process $Z=\{Z_t\}_{t\in \Z} $ is said to be $\mathcal{C}$-mixing see \cite{maume2006exponential} and \cite{hang2017bernstein}, if it satisfies
\[
\phi_{\mathcal{C}}(k) := \sup_{j\geq 1} \Big\{ \big| \E\big( Z f(Z_{j+k}) \big) - \E(Z)\E\big(f(Z_{j+k}) \big) \big|:  Z \in  L_1(\sigma_{1}^ j, P), \|Z\|_1 \leq 1, f \in \mathcal{C}(\mz), \|f\|_{\mathcal{C}} \leq 1 \Big\} \underset{k \rightarrow}{\longrightarrow}  0,
\]
where $\phi_{\mathcal{C}}(k)$ denotes the $C$-mixing coefficients of the process.
For more details on such dependence structure and for examples of $\mathcal{C}$-mixing processes, see for instance \cite{maume2006exponential}, \cite{hang2016learning}, \cite{hang2017bernstein}.   
Let us note for example that, $\phi$-mixing processes are $\mathcal{C}$-mixing.
\begin{enumerate}
\item \textbf{Geometrically $\mathcal{C}$-mixing processes.}
The stochastic process $Z$ is said to be geometrically $\mathcal{C}$-mixing, if the $\mathcal{C}$-mixing coefficient satisfies
\begin{equation}\label{equa_phi_mixing_coef}
\phi_{\mathcal{C}}(k) \leq c\exp\big(-bk^{\varrho}\big),  ~~ \text{ for all } k\ge 1,   
\end{equation}
for some constants $c > 0$, $b > 0$, and $\varrho > 0$.

\medskip

\noindent
From \cite{hang2016learning}, one can find $n_0 \geq 1$ such that, the assumption \textbf{(A4)} is satisfied for $n \geq n_0$ with $C = 2$, $c_{\gamma} = 8$, $c_A = 8/3$, and $\varphi(n) = n/\big(\log n\big)^ {2/\varrho}$.
Hence, for $\kappa \geq 1$ given in \textbf{(A3)} and for geometrically $\mathcal{C}$-mixing processes:
\begin{itemize}
\item The convergence rate of the expected excess risk of the NPDNN and the SPDNN estimators is of order $\mathcal{O} \Big(n^{-\frac{\kappa s}{\kappa s + d}} \big(\log n \big)^ {3 + 2/\varrho} \Big)$ on the class of H\"older smooth functions.
\item This rate is of order  $\mathcal{O} \Big( \big(\phi_{n}^{\kappa/2} \lor \phi_{n} \big)(\log n)^ {3 + (2/\varrho) \lor \kappa}\Big)$ when we deal with the class of composition H\"older functions. 
\end{itemize}
When $\kappa=2$, these rates are minimax optimal (up to a logarithmic factor).
\item \textbf{Polynomially $\mathcal{C}$-mixing processes.}
The process $Z$ is said to be polynomially $\mathcal{C}$-mixing, if the mixing coefficient $\phi_{\mathcal{C}}(k)$ satisfies
\begin{equation}\label{equa_poly_C-mixing_coef}
\phi_{\mathcal{C}}(k) \leq c\cdot k^{-\varrho},  ~~ \text{ for all }  k \ge 1,   
\end{equation}
for some constants $c > 0$, and $\varrho > 0$.
\medskip

\noindent
 When $\varrho >2$, there exist $n_0 > 0$ such that Assumption \textbf{(A4)} holds for $n \geq n_0$ with $C = 2$, $c_{\gamma} = 8$, $c_A = 8/3$, and $\varphi(n) = n^{(\varrho - 2)/(\varrho + 1)}$, see \cite{hang2016learning}. 
 In this case, for $\kappa \geq 1$ given in \textbf{(A3)} and for polynomially $\mathcal{C}$-mixing processes, the application of the bounds above gives:
\begin{itemize}
\item The convergence rate of the expected excess risk of the NPDNN and the SPDNN estimators is of order $\mathcal{O} \Big(n^{-\frac{\varrho - 2}{\varrho + 1}\frac{\kappa s}{\kappa s + d}} \big(\log n \big)^3 \Big)$, when the target function belongs to the class of H\"older smooth functions.
\item This convergence rate is of order  $\mathcal{O} \Big( \big(\phi _{n} ^{\frac{\kappa (\varrho - 2)}{2(\varrho + 1)}} \lor \phi _{n}^{\frac{\varrho-2}{\varrho + 1}} \big) (\log n)^3 \Big)$, on the class of composition H\"older functions. 
\end{itemize}
\end{enumerate}
\section{Nonparametric autoregression with exogenous covariate}\label{sect_np_reg}
Consider the nonparametric autoregression model given by
\begin{equation}\label{nonparametric_autoreg_model}
Y_t = f(Y_{t-1}, Y_{t-2}, \cdots, Y_{t-p}; \mx_{t-1}, \cdots, \mx_{t-q}) + \varepsilon_t,  ~ t \in \Z  
\end{equation}
where $f : \R^p \times \R^{q} \to \R$, (with $p, q \in \N$) is a measurable function, $(\mx_t)_{t \in \Z}$ is a process of exogenous covariates, and $(\varepsilon_t)_{t\in \Z}$ is a real-valued and centered sequence of i.i.d. random variables.
It is assumed an autoregressive-type structure on the covariates, that is,
\begin{equation}\label{cond_funct_g}
 \mx_t = g(\mx_{t-1},\cdots,\mx_{t-q}) + \eta_t  
 \end{equation}
where $(\eta_t)_{t \in \Z}$  is a real-valued and centered sequence of i.i.d. random variables, and $g : \R^q \to \R$ is a measurable function.
Set
\[   X_t = (Y_{t-1}, Y_{t-2}, \cdots, Y_{t-p}; \mx_{t-1}, \cdots, \mx_{t-q})^T , \xi_t = (\varepsilon_t, 0, \cdots ,0,\eta_t, 0, \cdots, 0)^T \in \R^{p+q} \text{ and }  u_t=(\varepsilon_t , \eta_t).\]
We have for all $t \in \Z$,
\begin{equation}\label{Markov_repre}
  X_{t+1} = \psi(X_t) + \xi_t, 
 \end{equation}
 where 
 $ \psi(y_1,\cdots,y_p; x_1,\cdots x_q) = \big( f(y_1,\cdots,y_p ; x_1,\cdots x_q), y_1,\cdots,y_{p-1}, g(x_1,\cdots, x_q), x_1,\cdots, x_{q-1}  \big)^T $ for all \\
$(y_1,\cdots,y_p; x_1,\cdots x_q) \in \R^{p+q}$. 

In the sequel, we need the following conditions:
\begin{itemize}
\item[\textbf{(B1)}:] The process $(u_t)_{t\in \Z}$ is i.i.d. and $u_0$ has a density with respect to the Lebesgue measure, which is positive and uniformly continuous on $\R$, and $\E[\|u_0\|^2] <\infty$.
\item[\textbf{(B2)}:] $|f(0;0)| < \infty$ and there exists non-negative constants $\alpha_{j,Y}(f)$ $j=1,\cdots,p$ and $\alpha_{j,\mx}(f)$, $\alpha_{j}(g)$ $j=1,\cdots,q$ such that the functions $f$ and $g$ given in (\ref{nonparametric_autoreg_model}) and (\ref{cond_funct_g}) satisfy for all $(y_1,\cdots,y_p), (y'_1,\cdots,y'_p) \in \R^p$, $(x_1,\cdots,x_q), (x'_1,\cdots,x'_q) \in \R^q$,
\[ |f(y_1,\cdots,y_p; x_1,\cdots,x_q) - f(y'_1,\cdots,y'_p; x'_1,\cdots,x'_q)| \leq \sum_{j=1}^p \alpha_{j,Y}(f) |y_j - y'_j| + \sum_{j=1}^q \alpha_{j,\mx}(f) |x_j - x'_j|,  \]
and
\[ |g( x_1,\cdots,x_q) - g(x'_1,\cdots,x'_q)| \leq \sum_{j=1}^q \alpha_{j}(g) |x_j - x'_j|. \]
Furthermore, 
$  \sum_{j=1}^p \alpha_{j,Y}(f) < 1$ and $\sum_{j=1}^q \alpha_{j}(g) < 1$.
\end{itemize}
\medskip

\noindent
 Assumption \textbf{(B2)} is a Lipschitz-type condition on the functions $f$ and $g$, and is classical in the study of the stability properties of such models, see for instance \cite{diop2022inference}.
Under \textbf{(B1)}, \textbf{(B2)} and from (\ref{Markov_repre}), there exists a stationary and ergodic geometrically $\alpha$-mixing solution $(Y_t, \mx_t)_{t \in \Z}$ of (\ref{nonparametric_autoreg_model}), see \cite{chen2000geometric} and \cite{doukhan1994mixing}. 

\medskip
Consider the nonparametric estimation of the function  $f^{*} = f\ind_{\mx}$ in (\ref{nonparametric_autoreg_model}) by the NPDNN and the SPDNN predictors. 
Set,
\[X_t = (Y_{t - 1}, \cdots, Y_{t - p}; \mathcal{X}_{t - 1}, \cdots, \mx_{t-q}).  \]
 A target predictor according to the squared loss function is given for all $x \in \R^{p +q} $ by
\begin{equation}
h^{*} (x) = \E[ Y_0 | X_0 =x] = f(x).
\end{equation}
 If the distribution of $\varepsilon_0$ is symmetric around 0, then a target function with respect to the $L_1$ loss and the Huber loss is given respectively by (see \cite{fan2024noise}, \cite{shen2021deep}),
\begin{equation*}
h^{*} (x) = \text{med} ( Y_0 | X_0 = x) = f(x) \text{ and } h^{*}(x) = f(x),
\end{equation*}
for all $x \in \R^{p +q}$, where $\text{med}(V)$ denotes the median of $V$ for any random variable $V$. 
In this case, with the Huber loss (where \textbf{(A3)} holds with $\kappa = 2$), the convergence rate of the expected excess risk of the NPDNN and SPDNN estimator are $\mathcal{O} \Big(n^{-\frac{2 s}{2s + d}} \big( \log (n) \big)^5 \Big)$ and  $\mathcal{O} \big(\phi_n (\log n)^{5}\big)$ respectively on the class of H\"older functions (with a smoothness parameter $s>0$) and composition H\"older functions defined in (\ref{equa_compo_structured_function}).

\section{Proofs of the main results}\label{prove}
\subsection{Proof of Theorem  \ref{thm2}} 
In this proof, we follow the same steps as in the proof of Theorem 3.1 in \cite{kengne2025robust}.
Recall the training set $ D_n = \{ Z_1 = (X_1, Y_1), \cdots, Z_n= (X_n, Y_n) \} $.
 Consider another sample  $ D'_n \coloneqq  \{ Z'_1 = (X'_1, Y'_1), \cdots, Z'_n= (X'_n, Y'_n) \} $ generated from $\{Z_t =(X_t, Y_t), t \in \Z \} $ and independent of $D_n$.
Let us denote for all predictor $h:\mathcal{X} \rightarrow \mathcal{Y}$ and $i=1,\cdots,n$, 
\begin{equation}\label{proof_def_g}
g(h(X_i), Y_i) := \ell(h(X_i), Y_i) - \ell(h^{*}(X_i), Y_i). 
\end{equation} 
Let $L_n, N_n, B_n, F_n, S_n >0$ fulfilling the conditions in Theorem \ref{thm2}.
Set:
\begin{equation}\label{proof_def_H_sigma_n}
 \mathcal{H}_{\sigma,n}  := \mathcal{H}_{\sigma}(L_n, N_n, B_n, F_n, S_n) \text{ and }  h_{\mathcal{H}_{\sigma, n}} := \underset{h \in \mathcal{H}_{\sigma, n}}{\argmin} R(h).
\end{equation} 
By using the equation $\E[R(\widehat{h}_{n,NP}) - R(h^{*})] = \E_{D_n} \Big[\E_{D_n'} \Big[\frac{1}{n}\sum_{i= 1}^n g(\widehat{h}_{n,NP}, Z_i') \Big] \Big]$ and the inequality $\widehat{R}_n(\widehat{h}_{n,NP}) \leq \widehat{R}_n(h_{\mathcal{H}_{\sigma, n}})$, one can easily get
\begin{equation}\label{Exp_excess_risk_v2}
  \E[R(\widehat{h}_{n,NP}) - R(h^{*})] \leq  \E_{D_n} \Big[\dfrac{1}{n}\sum_{i= 1}^n \{ -2g(\widehat{h}_{n,NP}, Z_i)+ \E_{D_n'} g(\widehat{h}_{n,NP}, Z_i') \} \Big] + 2 [R(h_{\mathcal{H}_{\sigma, n}}) - R(h^{*})]. 
\end{equation}

\medskip

\noindent
\textbf{Step 1 : Bounding the first term in the right-hand side of (\ref{Exp_excess_risk_v2})}.
For any $h \in \mathcal{H}_{\sigma, n}$ and $ \epsilon >0$, set
\[G(h, Z_i) :=  \E_{D_n'} g(h, Z_i') -2g(h, Z_i),
\]
and 
 \begin{equation}\label{proof_def_m_n}
 m := \mathcal{N}\Big( \mathcal{H}_{\sigma, n}, \epsilon  \Big).
\end{equation} 
This covering number is finite, $\|h\|_\infty \leq F_n < \infty$ for all $h \in \mathcal{H}_{\sigma, n}$. 
Let $h_1,\cdots,h_m \in \mathcal{H}_{\sigma,n}$ such that,
\begin{equation}\label{proof_H_sigma_n_subset_ball}
\mathcal{H}_{\sigma,n} \subset \bigcup_{j=1}^m B(h_j, \epsilon),
\end{equation}
where $B(h_j, \epsilon)$ is the ball of radius $\epsilon$, centered at $h_j$.
So, there exists a random index $j^* \in \{1,\cdots, m\}$ such that $ \|\widehat{h}_{n,NP} - h_{j^*}\|_{\infty} \leq \epsilon$.
From the Lipschitz continuity of the loss function, we have for all $i=1,\cdots,n$, 
\begin{align*}
\nonumber |g(\widehat{h}_{n,NP}), Z_i) - g(h_{j^*}, Z_i)| & = |\ell(\widehat{h}_{n,NP}(X_i), Y_i) - \ell(h_{j^*}(X_i), Y_i)|  \leq \mathcal{K}_{\ell} |\widehat{h}_{n,NP}(X_i) - h_{j^*}(X_i)| \leq \mathcal{K}_{\ell} \epsilon.
\end{align*}
 One can easily get,
\begin{equation*}
  \dfrac{1}{n}\sum_{i=1}^n\E_{D_n}[g(\widehat{h}_{n,NP}, Z_i)]   \leq  \dfrac{1}{n}\sum_{i=1}^n\E_{D_n}[g(h_{j^*}, Z_i)] + \mathcal{K}_{\ell} \epsilon,
\end{equation*}
and
\begin{align}\label{euqua_bound_estimation_error}
 \nonumber \E_{D_n} \Big[\dfrac{1}{n}\sum_{i= 1}^n \{ -2g(\widehat{h}_{n,NP}, Z_i)+ \E_{D_n'} g(\widehat{h}_{n,NP}, Z_i') \} \Big] &= \E_{D_n}\Big[\dfrac{1}{n}\sum_{i=1}^n G(\widehat{h}_{n,NP}, Z_i) \Big] \\
 \nonumber &\leq \E_{D_n}\Big[\dfrac{1}{n}\sum_{i=1}^n G(h_{j^*}, Z_i) \Big] + \E_{D_n}\Big[\dfrac{1}{n}\sum_{i=1}^n |G(\widehat{h}_{n,NP}, Z_i) - G(h_{j^*}, Z_i)|   \Big] \\
 &\leq \E_{D_n}\Big[\dfrac{1}{n}\sum_{i=1}^n G(h_{j^*}, Z_i) \Big] + 3 \mathcal{K}_{\ell} \epsilon.   
\end{align}
From the Proposition 1 in \cite{ohn2019smooth}, we have,
\begin{equation}\label{proof_ing_covering_num}
\mathcal{N} \Big (\mathcal{H}_{\sigma, n}, \epsilon \Big)  \leq  \exp \Bigg (  2 L_n (S_n + 1) \log \left(\frac{1}{\epsilon} C_{\sigma} L_n (N_n + 1)(B_n \lor 
 1) \right) \Bigg),
\end{equation}
where $B_n \lor 1 = \max(B_n, 1)$.
The dependence of $m$ on $n$ is omitted to simplify notation.

\medskip
\noindent
Recall that $\|h^{*}\|_{\infty} \leq \mathcal{K} $ and for all $x \in \mathcal{X}$,
\[h^{*}(x) = \underset{h(x): \|h\|_{\infty} \leq \mathcal{K} }{\argmin} \E[\ell(h(X), Y)|X=x].
\]
Let $ h \in \mathcal{H}_{\sigma, n}$.
Since $\|h\|_\infty \leq F_n$, according to Assumption \textbf{(A2)}, we have,
\begin{align}\label{proof_g_beta_K_ell_beta}
  |g(h, Z_i)| & = |\ell(h(X_i), Y_i) - \ell(h^{*}(X_i), Y_i)| \leq \mk_{\ell} | h(X_i) - h^{*}(X_i)| \leq \mk_{\ell}(F + \mk).
\end{align}
Hence,
\begin{align}\label{equa_bound_of_the_variance}
  \var(g(h, Z_i))  \leq \E[g(h, Z_i)^2]  \leq \E[|g
  (h, Z_i)||g(h, Z_i)|]
  \leq  \mk_{\ell}(F + \mk) \E[g(h, Z_i)]. 
\end{align}
Also, we have from (\ref{proof_g_beta_K_ell_beta}),  $|\E_{D_n'}[g(h_{j^*}, Z_i')] - g(h_{j^*}, Z_i)| \leq 2\mk_{\ell}(F + \mk)$. 
Let $ \varepsilon > 0$.
According to (\ref{equa_bound_of_the_variance}) and the generalized Bernstein type inequality given in (\ref{G_beernstein_expo_inequa}), we have for all $j=1,\cdots,m$ (defined in (\ref{proof_def_m_n})),
\begin{align}\label{term_bound_excess_risk}
\nonumber P \Bigg\{\dfrac{1}{n} \sum_{i=1}^n G(h_{j}, Z_i) > \varepsilon \Bigg\}  & = P \Bigg\{\dfrac{1}{n} \sum_{i=1}^n\E_{D_n'}[g(h_{j}, Z_i')] - \dfrac{2}{n} \sum_{i = 1}^n g(h_{j}, Z_i) > \varepsilon \Bigg\}
\\
\nonumber & =  P \Big\{\E_{D_n'}[g(h_{j}, Z_1')] - \dfrac{2}{n} \sum_{i = 1}^n g(h_{j}, Z_i) > \varepsilon \Big\}.
\\
\nonumber & =  P \Big\{2\E_{D_n'}[g(h_{j}, Z_1')] - \dfrac{2}{n} \sum_{i = 1}^n g(h_{j}, Z_i) > \varepsilon + \E_{D_n'}[g(h_{j}, Z_1')] \Big\}
\\
\nonumber & =  P \Big\{\E_{D_n'}[g(h_{j}, Z_1')] - \dfrac{1}{n} \sum_{i = 1}^n g(h_{j}, Z_i) > \dfrac{\varepsilon}{2} + \dfrac{1}{2} \E_{D_n'}[g(h_{j}, Z_1')]\Big\}
\\
\nonumber & \leq  P \Bigg\{\E_{D_n'}[g(h_{j}, Z_1')] - \dfrac{1}{n} \sum_{i = 1}^n g(h_{j}, Z_i) > \dfrac{\varepsilon}{2} + \dfrac{\E_{D_n'}[g(h_{j}, Z_1')]}{2} \Bigg\}
\\
\nonumber & \leq  P \Bigg\{\dfrac{1}{n} \sum_{i = 1}^n \Big(\E[g(h_{j}, Z_i)] - g(h_{j}, Z_i)\Big) > \dfrac{\varepsilon}{2} + \dfrac{\E[g(h_{j}, Z_1)]}{2} \Bigg\}
\\
\nonumber & \leq C \exp\Bigg[-\dfrac{\Big(\dfrac{\varepsilon}{2} + \dfrac{\E[g(h_{j}, Z_1)]}{2} \Big)^2 \varphi(n)}{c_\gamma \mk_{\ell}(F + \mk) \E[g(h, Z_i)] +  c_A\Big(\dfrac{\varepsilon}{2} + \dfrac{\E[g(h_{j}, Z_1)}{2} \Big) \big(2\mk_{\ell}(F + \mk) \big)} \Bigg]
\\
\nonumber & \leq C\exp\Bigg[-\dfrac{(1/4)\Big(\varepsilon + \E[g(h_{j}, Z_1)] \Big)^2 \varphi(n)}{c_\gamma \mk_{\ell}(F + \mk) \big(\varepsilon + \E[g(h, Z_i)] \big) +  c_A\Big(\dfrac{\varepsilon}{2} + \dfrac{\E[g(h_{j}, Z_1)}{2} \Big) \big(2\mk_{\ell}(F + \mk) \big)}  \Bigg]
\\
\nonumber & \leq C\exp\Bigg[-\dfrac{(1/4)\Big(\varepsilon + \E[g(h_{j}, Z_1)] \Big) \varphi(n)}{c_\gamma \mk_{\ell}(F + \mk)  +  c_A\mk_{\ell}(F + \mk)}  \Bigg]
 \leq C\exp\Bigg[-\dfrac{\Big(\varepsilon + \E[g(h_{j}, Z_1)] \Big) \varphi(n)}{(1/4)(c_\gamma  +  c_A)\mk_{\ell}(F + \mk)}  \Bigg]
\\
& \leq C\exp\Bigg[-\dfrac{\varepsilon \varphi(n)}{(1/4)(c_\gamma  +  c_A)\mk_{\ell}(F + \mk)}  \Bigg],
\end{align}
where $c_\gamma$ and $c_A$ are positive constants given in (\ref{G_beernstein_expo_inequa}). In the sequel, we set $\widetilde{C} := \widetilde{C}(\mk_{\ell}, F_n, \mk) = (1/4)(c_\gamma + c_A)$.
Set, $\epsilon = \frac{1}{\varphi(n)}$ in (\ref{proof_def_m_n}) and let $m = \mathcal{N}\Big( \mathcal{H}_{\sigma, n},  1/ \varphi(n) \Big)$. We get, in addition to (\ref{proof_ing_covering_num}),
\begin{align*}
\nonumber P \Big\{\dfrac{1}{n} \sum_{i = 1}^n G(h_{j^*}, Z_i) > \varepsilon \Big\}  & \leq   P \Bigg\{ \bigcup_{j=1}^m \Big( \dfrac{1}{n} \sum_{i = 1}^n G(h_{j}, Z_i) > \varepsilon \Big) \Bigg\} \leq \sum_{j=1}^m  P \Big\{\dfrac{1}{n} \sum_{i = 1}^n G(h_{j}, Z_i) > \varepsilon \Big\}  
\\
\nonumber & \leq \mathcal{N}\Big( \mathcal{H}_{\sigma, n}, \frac{1}{\varphi(n)}\Big)\cdot C\exp\Bigg[-\dfrac{\varepsilon \varphi(n)}{\widetilde{C}\mk_{\ell}(F + \mk)}  \Bigg]
\\
& \leq C\exp \Bigg[2L_n(S_n + 1)\log\Big(\varphi(n) C_{\sigma}L_n(N_n + 1)(B_n \lor 1)\Big) -\dfrac{\varepsilon \varphi(n)}{\widetilde{C}\mk_{\ell}(F + \mk)} \Bigg]. 
\end{align*}
For any $\alpha_n > 0$, we have
\begin{align*} %\label{equa1_bound_G_beta_n}
\nonumber & \E\Big[\dfrac{1}{n}\sum_{i = 1}^nG(h_{j^*}, Z_i) \Big]  \leq \alpha_n + \displaystyle\int_{\alpha_n}^{\infty} P\Big\{\dfrac{1}{n}\sum_{i = 1}^nG(h_{j}, Z_i) > \varepsilon\Big\} d\varepsilon \\
\nonumber & \leq \alpha_n + C\exp \Bigg[2L_n(S_n + 1)\log\Big(\varphi(n) C_{\sigma}L_n(N_n + 1)(B_n \lor 1)\Big) \Bigg] \times \displaystyle\int_{\alpha_n}^{\infty} \exp \Bigg( -\dfrac{\varepsilon \varphi(n)}{\widetilde{C}\mk_{\ell}(F + \mk)} \Bigg) d\varepsilon 
\\
\nonumber & \leq \alpha_n +
\frac{\widetilde{C}\mk_{\ell}(F + \mk)}{\varphi(n)}\exp\Bigg[2L_n(S_n + 1)\log \Big(\varphi(n) C_{\sigma}L_n(N_n + 1)(B_n \lor 1) \Big)  \Bigg] \exp \Bigg[-\dfrac{\alpha_n \varphi(n)}{\widetilde{C}\mk_{\ell}(F + \mk)} \Bigg]
\\
 & \leq \alpha_n +
\frac{\widetilde{C}\mk_{\ell}(F + \mk)}{\varphi(n)}\exp\Bigg[2L_n(S_n + 1)\log \Big(\varphi(n) C_{\sigma}L_n(N_n + 1)(B_n \lor 1) \Big) -\dfrac{\alpha_n \varphi(n)}{\widetilde{C}\mk_{\ell}(F + \mk)} \Bigg].
\end{align*}
Hence, with 
\begin{equation}\label{equa_alpha_n}
\alpha_n \coloneqq \dfrac{ \big(\log \big(\varphi(n) \big)  \big)^{\nu}}{ \big(\varphi(n) \big)^{\frac{\kappa s}{\kappa s + d} } },  
\end{equation} 
for some $\nu >3$, we get
\begin{align}\label{proof_E_G_beta_log_nu_n}
\nonumber &\E\Big[ \dfrac{1}{n}\sum_{i = 1}^nG(h_{j^*}, Z_i) \Big] \\
& \leq \dfrac{ \big(\log \big(\varphi(n) \big)  \big)^{\nu}}{ \big(\varphi(n) \big)^{\frac{\kappa s}{\kappa s + d} } } + \frac{\widetilde{C}\mk_{\ell}(F + \mk)}{\varphi(n)}\exp\Bigg[2L_n(S_n + 1)\log \Big(\varphi(n) C_{\sigma}L_n(N_n + 1)(B_n \lor 1) \Big) -\dfrac{ \big(\varphi(n) \big)^{\frac{d}{\kappa s + d}} \big(\log \big(\varphi(n) \big)  \big)^{\nu}}{\widetilde{C}\mk_{\ell}(F + \mk)} \Bigg].
\end{align}
Since
$
L_n = \dfrac{s L_0}{\kappa s + d} \log \big(\varphi(n)\big), N_n = N_0 \big(\varphi(n)\big)^{\frac{d}{\kappa s + d}}, S_n = \frac{s S_0}{\kappa s + d} \big(\varphi(n)\big)^{\frac{d}{\kappa s + d}}\log\big(n^{(\alpha)} \big), 
$
$ B_n = B_0 \big(\varphi(n) \big)^{\frac{4(d + s)}{\kappa s + d}}$,
with $L_0, N_0, B_0, S_0 >0$,
 it holds from (\ref{proof_E_G_beta_log_nu_n}) that,
\begin{align}\label{equa2_bound_G_beta_n}
\nonumber & \E\Big[\dfrac{1}{n}\sum_{i = 1}^nG(h_{j^*}, Z_i) \Big] \\
\nonumber &  \leq \dfrac{ \big(\log \big(\varphi(n)\big)  \big)^{\nu}}{ \big(\varphi(n)\big)^{\frac{\kappa s}{\kappa s + d}} } + \frac{\widetilde{C}\mk_{\ell}(F + \mk)}{\varphi(n)} \exp\Bigg[\dfrac{ \big(\varphi(n) \big)^{\frac{d}{\kappa s + d}} \big(\log \big(\varphi(n) \big)  \big)^{\nu}}{\widetilde{C}\mk_{\ell}(F + \mk)} \\
\nonumber & \hspace{8cm} \times   \Bigg(\dfrac{2\widetilde{C}\mk_{\ell}(F + \mk) \dfrac{s L_0}{\kappa s + d} \log\big(\varphi(n)\big) \Big(\frac{s S_0}{\kappa s + d}\log\big(\varphi(n)\big) + 1 \Big)}{ \big(\log\big(\varphi(n)\big) \big)^{\nu}} \\
& \hspace{2cm}   \times\log\Big(\varphi(n) C_{\sigma} \dfrac{s L_0}{\kappa s + d} \log\big(\varphi(n)\big) \big(N_0 \big(\varphi(n) \big)^{\frac{d}{\kappa s + d}} + 1 \big) \big(B_0 \big(\varphi(n) \big)^{\frac{4(d + s)}{\kappa s + d}} \lor 1 \big) \Big)
- 1 \Bigg) \Bigg].
\end{align}

\medskip

Since $\nu >3$, we have
\begin{align*}
&\dfrac{2\widetilde{C}\mk_{\ell}(F + \mk) \dfrac{s L_0}{\kappa s + d} \log\big(\varphi(n)\big) \Big(\frac{s S_0}{\kappa s + d}\log\big(\varphi(n)\big) + 1 \Big)}{ \big(\log\big(\varphi(n)\big) \big)^{\nu}} \\
& \hspace{2cm}   \times\log\Big(\varphi(n) C_{\sigma} \dfrac{s L_0}{\kappa s + d} \log\big(\varphi(n)\big) \big(N_0 \big(\varphi(n) \big)^{\frac{d}{\kappa s + d}} + 1 \big) \big(B_0 \big(\varphi(n) \big)^{\frac{4(d + s)}{\kappa s + d}} \lor 1 \big) \Big) \limiten 0.
\end{align*}
So, one can find $n_0 = n_0(\mk_\ell, \mk, F, \kappa, L_0, N_0, B_0, S_0, C_\sigma, s, d,)$ such that, for $n \ge n_0$, we get from (\ref{equa2_bound_G_beta_n}),
\begin{align*} %\label{equa2_bound_G_beta_n_j_star}
\nonumber \E\Big[\dfrac{1}{n}\sum_{i = 1}^nG(h_{j^*}, Z_i) \Big] & \leq \dfrac{ \big(\log \big(\varphi(n)\big)  \big)^{\nu}}{ \big(\varphi(n)\big)^{\frac{\kappa s}{\kappa s + d}} } + \frac{\widetilde{C}\mk_{\ell}(F + \mk)}{\varphi(n)} \exp\Bigg(- \dfrac{ \big(\varphi(n) \big)^{\frac{d}{\kappa s + d}} \big(\log \big(\varphi(n) \big)  \big)^{\nu}}{2\widetilde{C}\mk_{\ell}(F + \mk)} \Bigg) \\
&\leq \dfrac{ \big(\log \big(\varphi(n)\big)  \big)^{\nu}}{ \big(\varphi(n)\big)^{\frac{\kappa s}{\kappa s + d}} } + \frac{\widetilde{C}\mk_{\ell}(F + \mk)}{\varphi(n)}.
\end{align*}
Hence, in addition to (\ref{euqua_bound_estimation_error}) and since $\epsilon = 1/\varphi(n)$, we have
\begin{equation} \label{equa_bound_stochastic_error_v2}
 \E_{D_n}  \Big[\dfrac{1}{n}\sum_{i= 1}^n \{ -2g(\widehat{h}_{n,NP}, Z_i)+ \E_{D_n'} g(\widehat{h}_{n,NP}, Z_i') \} \Big]  \leq   \dfrac{ \big(\log \big(\varphi(n)\big)  \big)^{\nu}}{ \big(\varphi(n)\big)^{\frac{\kappa s}{\kappa s + d}} } + \frac{\widetilde{C}\mk_{\ell}(F + \mk)}{\varphi(n)} + \frac{3\mathcal{K}_{\ell}}{\varphi(n)}.    
\end{equation} 

\medskip

\noindent
\textbf{Step 2 : Bounding the second term in the right-hand side of (\ref{Exp_excess_risk_v2})}. 

\medskip

\noindent
Set $\epsilon_n = \dfrac{ 1}{ \big(\varphi(n)\big)^{\frac{s}{\kappa s + d}} } $.
According to Theorem 3.2 in \cite{kengne2025excess}, there exists a positive constants $L_0, N_0, B_0, S_0 >0$ such that with 
$L_n = \dfrac{s L_0}{\kappa s + d} \log\big(\varphi(n)\big), N_n = N_0 \big(\varphi(n)\big)^{\frac{d}{ \kappa s + d}}, S_n = \frac{s S_0}{\kappa s + d} \big(\varphi(n)\big)^{\frac{d}{\kappa s + d}}\log\big(\varphi(n)\big), 
$
$ B_n = B_0 \big(\varphi(n)\big)^{\frac{4(d + s)}{\kappa s + d}}$, there is a neural network $h_{n} \in \mathcal{H}_{\sigma, n} = \mathcal{H}_{\sigma}(L_n, N_n, B_n, F_n, S_n)$ satisfying,
\begin{equation}\label{cond_excess_risk_bound_np}
\| h_{n} - h^{*} \|_{\infty, \mathcal{X}} \leq \epsilon_n= \dfrac{ 1}{ \big( \varphi(n)\big)^{\frac{s}{\kappa s + d} } }.
\end{equation}
Set
\[   \widetilde{\mathcal{H}}_{\sigma, n} := \{h \in  \mathcal{H}_{\sigma,n}, ~ \|h - h^*\|_{\infty, \mx}  \leq \epsilon_n\}  .\]
For $\varphi(n) \geq \varepsilon_0^{-( \kappa + d/s )}$, we get from (\ref{cond_excess_risk_bound_np}) and \textbf{(A3)} that,
\begin{align}\label{approxi_error2}
\nonumber R(h_{\mathcal{H}_{\sigma, n}}) - R(h^{*}) & = \underset{h\in \mathcal{H}_{\sigma, n}}{\inf} R(h) - R(h^{*}) \leq \underset{h\in  \widetilde{\mathcal{H}}_{\sigma, n}}{\inf} \Big(R(h) - R(h^{*}) \Big)
\\
& \leq \mathcal{K}_{0} \underset{h\in  \widetilde{\mathcal{H}}_{\sigma, n}}{\inf}  \| h - h^*\|^\kappa_{\kappa, P_{X_0}} \leq \dfrac{ \mk_0}{ \big( \varphi(n) \big)^{\frac{\kappa s}{\kappa s + d}} }.
\end{align}

\medskip

\noindent
\medskip

\noindent
\textbf{Step 3 : Expected excess risk bound}.
\medskip

\noindent
According to (\ref{Exp_excess_risk_v2}), (\ref{equa_bound_stochastic_error_v2}) and (\ref{approxi_error2}), we have for $n \geq n_0$,
\begin{align*}
\E[R(\widehat{h}_{n,NP}) - R(h^{*})]  
& \leq  \dfrac{ \big(\log \big(\varphi(n)\big)  \big)^{\nu}}{ \big(\varphi(n)\big)^{\frac{\kappa s}{\kappa s + d}} } + \frac{\widetilde{C}\mk_{\ell}(F + \mk)}{\varphi(n)} + \frac{3\mathcal{K}_{\ell}}{\varphi(n)} + \dfrac{2\mk_0}{ \big(\varphi(n)\big)^{\frac{\kappa s}{\kappa s + d}} } \\
& \leq \dfrac{ \big(\log \big(\varphi(n)\big)  \big)^{\nu} + 2\mk_0}{ \big(\varphi(n)\big)^{\frac{\kappa s}{\kappa s + d}} } + \frac{\widetilde{C}\mk_{\ell}(F + \mk) + 3 \mk_0  }{  \varphi(n) } \lesssim \dfrac{ \big(\log \big(\varphi(n)\big)  \big)^{\nu}}{ \big(\varphi(n)\big)^{\frac{\kappa s}{\kappa s + d}} } ,
\end{align*}
for all $\nu >3$.

\qed

\subsection{Proof of Theorem \ref{thm3}}
%o
Recall the following decomposition of the expected excess risk in (\ref{Exp_excess_risk_v2}), see Proof of Theorem \ref{thm2} above:
\begin{align}\label{Exp_excess_risk_v3}
\E[R(\widehat{h}_{n,NP}) - R(h^{*})] \leq  \E_{D_n} \Big[\dfrac{1}{n}\sum_{i= 1}^n \{ -2g(\widehat{h}_{n,NP}, Z_i)+ \E_{D_n'} g(\widehat{h}_{n,NP}, Z_i') \} \Big] + 2 [R(h_{\mathcal{H}_{\sigma, n}}) - R(h^{*})], 
\end{align}
with 
\begin{equation*} 
 \mathcal{H}_{\sigma,n}  := \mathcal{H}_{\sigma}(L_n, N_n, B_n, F_n, S_n) \text{ and }  h_{\mathcal{H}_{\sigma, n}} := \underset{h \in \mathcal{H}_{\sigma, n}}{\argmin} R(h),
\end{equation*} 
where $L_n, N_n, B_n, F_n, S_n >0$ fulfilling the conditions of this theorem.
%
% %
Using similar arguments as in \textbf{Step 1} in the proof of Theorem \ref{thm2}, and by setting $ \alpha_n \coloneqq \big(\log \big(\varphi(n) \big)  \big)^{\nu}\phi_{n, \varphi}$ as in (\ref{equa_alpha_n}), for some $\nu >3$, and $\phi_{n, \varphi}$ defined in (\ref{def_phi_nalpha}), we obtain the following bound for the first term in the right-hand side of (\ref{Exp_excess_risk_v3}), 
\begin{equation} \label{equa_bound_stochastic_error_v3}
 \E_{D_n}  \Big[\dfrac{1}{n}\sum_{i= 1}^n \{ -2g(\widehat{h}_{n,NP}, Z_i)+ \E_{D_n'} g(\widehat{h}_n, Z_i') \} \Big]  \lesssim  \big(\log \big(\varphi(n) \big)  \big)^{\nu}\phi_{n, \varphi} .    
\end{equation} 
Let us derive a bound of the second term in the right-hand side of (\ref{Exp_excess_risk_v3}), on the class of composition H\"older functions $\mathcal{G}(q, \bold{d}, \bold{t}, \boldsymbol{\beta}, \mathcal{K})$.
According to the proof of Theorem 1 in \cite{schmidt2020nonparametric}, there exists a constant $C_0>0$ independent of $n$ such that, 
\begin{equation}\label{cond_excess_risk_bound_corol01_v1_2}
\underset{h^{*}\in \mathcal{G}(q, \bold{d}, \bold{t}, \boldsymbol{\beta}, A)}{\sup}~\underset{h \in \mathcal{H}_{\sigma, n}}{\inf}\|h - h^{*}\|_{\infty, \mathcal{X}}^2 \leq C_0 \underset{0\leq i \leq q}{\max}\big(\varphi(n) \big)^{ \frac{-2\beta_i^* }{2\beta_i^* + t_i}   } = C_0 \phi_{n, \varphi}.  
\end{equation}
Set
\[  \widetilde{\mathcal{H} }_{\sigma, n} := \{h \in  \mathcal{H}_{\sigma,n}, ~ \|h - h^*\|_{\infty, \mx}  \leq \sqrt{\phi_{n, \varphi}}\}.
\]
Thus, under \textbf{(A3)} and for sufficiently large $n$ such that $\sqrt{C _0 \phi_{n, \varphi} } \leq \varepsilon _0$, from (\ref{cond_excess_risk_bound_corol01_v1_2}),
for all $h^{*} \in \mathcal{G}(q, \bold{d}, \bold{t}, \boldsymbol{\beta}, A)$, we get,
\begin{align}\label{excess_risk_bound_compo_holder_funct}
\nonumber R(h_{\mathcal{H}_{\sigma, n}}) - R(h^{*}) & =  \underset{h \in \mathcal{H}_{\sigma, n} }{\inf} \Big[R(h) - R(h^{*})\Big]
\\
\nonumber & \leq \mk_\ell \underset{h \in \mathcal{H}_{\sigma,n} }{\inf} \Big[\E\big[|h(X_0)  - h^ *(X_0)|^\kappa\big]\Big]   \\    
\nonumber & \leq  \mk_\ell \underset{h \in \widetilde{\mathcal{H}}_{\sigma,n} }{\inf} \Big[\|h- h^ {*} \|_{\infty, \mathcal{X}} ^{\kappa} \\
 & \lesssim \phi_{n, \varphi} ^{\kappa /2}.
\end{align}
So, from (\ref{Exp_excess_risk_v3}), (\ref{equa_bound_stochastic_error_v3}), (\ref{excess_risk_bound_compo_holder_funct}), we have for all $h^{*} \in \mathcal{G}(q, \bold{d}, \bold{t}, \boldsymbol{\beta}, A)$ and $\nu >3$,
\begin{align*}
\E[R(\widehat{h}_{n,NP}) - R(h^{*})] \leq \big(\log \big(\varphi(n) \big)  \big)^{\nu}\phi_{n, \varphi} + \phi_{n, \varphi} ^{\kappa /2} 
\lesssim \Big(\phi _{n, \varphi} ^{\kappa/2} \lor \phi _{n, \varphi}\Big) \Big(\log \big( \varphi(n) \big)  \Big)^{\nu}. 
\end{align*}
% %

\qed

\subsection{Proof of Theorem \ref{thm1}}

For this proof, we need the following lemma.
\begin{lem}\label{Lemma_bernstein_inequa}
Let $(U_1,\cdots,U_n)$ be a trajectory of a stationary and ergodic $\mz$-valued process $\{ Z_t = (X_t, Y_t), t \in \Z \}$ that satisfies \textbf{(A4)}.  
Let $\mathcal{G}$ be a set of real-valued functions on $\mathcal{Z}$ such that, for some $C_1, C_2 \ge 0$, $|g(U_1) - \E[g(U_1)]| \leq C_1$ $a.s.$ and $\E[g(U_1)] ^2 \leq C_2 \E[g(U_1)]$ for any $g \in \mathcal{G}$. Then for all $\varepsilon > 0$ and $0 < u \leq 1$, it holds that,
\begin{equation}
P\Big\{ \sup_{g \in \mathcal{G}} \dfrac{\E[g(U_1)] - (1/n) \sum_{i = 1} ^n g(U_i)}{\sqrt{ E[g(U_1)] + \varepsilon }} \ge 4u \sqrt{\varepsilon} \Big\} \lesssim \mathcal{N}(\mathcal{G}, u\varepsilon) \cdot \exp \Bigg( -\dfrac{u^2 \varepsilon \varphi(n) }{c_{\gamma}C_2 +  c_A C_1 } \Bigg),
\end{equation}
where the function $\varphi$ and the positive constants  $c_{\gamma}, c_A$ are given in \textbf{(A4)}.
\end{lem}

\medskip

\noindent
\textbf{Proof of Lemma \ref{Lemma_bernstein_inequa} }
Under Assumption \textbf{(A4)}, this lemma can be easily  obtained in the same way as in the proof of Lemma 6.1 in \cite{kengne2025deep}.
\qed

\medskip

\medskip

\noindent
\textbf{Proof of Theorem \ref{thm1}}
Consider the following decomposition of the expected excess risk: 
\begin{equation}\label{proof_excess_risk_v1}
 \E[R(\widehat{h}_{n,SP}) - R(h^{*}) ]  \coloneqq \E[B_{1, n}] +  \E[B_{2, n}],
 \end{equation}
where,
\[
\begin{array}{llll}
  B_{1, n}   & =   [R(\widehat{h}_{n,SP}) - R(h^{*})]  - 2 [ \widehat{R}_n (\widehat{h}_{n,SP})  - \widehat{R}_n (h^{*}) ] - 2 J_n(\widehat{h}_{n,SP}); 
  \\
 B_{2, n}   &  =   2[\widehat{R}_n (\widehat{h}_{n,SP}) - \widehat{R}_n (h^{*})] + 2 J_n(\widehat{h}_{n,SP}).
\end{array}
\]
 Following the same steps as in the proof of the Theorem 4.1 in \cite{kengne2025deep} with $\varphi(n)$ instead $n^{(\alpha)}$ and by using Lemma \ref{Lemma_bernstein_inequa} above, one can find $n_0 >0$ such that, for $n \geq n_0$ we get,
\[   \E [B_{1, n}] \leq \int_ 0^ {\infty} P(B_ {1, n} > \rho) d\rho  \lesssim \dfrac{1}{\varphi(n)} ~  \text{ and } ~ 
\E[B_{2, n}] \leq 2 \underset{h \in \mathcal{H}_{\sigma}(L_n, N_n, B_n, F)}{\inf}\big[ R(h) - R(h^*) + J_n(h) \big] + \dfrac{1}{\varphi(n)}.   \]
Which gives the result of the theorem. 
\qed


\begin{thebibliography}{10}

\bibitem{alquier2025minimax}
{\sc Alquier, P., and Kengne, W.}
\newblock Minimax optimality of deep neural networks on dependent data via
  pac-bayes bounds.
\newblock {\em Electronic Journal of Statistics 19}, 2 (2025), 5895--5924.

\bibitem{bauer2019deep}
{\sc Bauer, B., and Kohler, M.}
\newblock On deep learning as a remedy for the curse of dimensionality in
  nonparametric regression.
\newblock {\em The Annals of Statistics 47}, 4 (2019), 2261--2285.

\bibitem{chen2000geometric}
{\sc Chen, M., and Chen, G.}
\newblock Geometric ergodicity of nonlinear autoregressive models with changing
  conditional variances.
\newblock {\em Canadian Journal of Statistics 28}, 3 (2000), 605--614.

\bibitem{dicker2013variable}
{\sc Dicker, L., Huang, B., and Lin, X.}
\newblock Variable selection and estimation with the seamless-l 0 penalty.
\newblock {\em Statistica Sinica\/} (2013), 929--962.

\bibitem{diop2022inference}
{\sc Diop, M.~L., and Kengne, W.}
\newblock Inference and model selection in general causal time series with
  exogenous covariates.
\newblock {\em Electronic Journal of Statistics 16}, 1 (2022), 116--157.

\bibitem{doukhan1994mixing}
{\sc Doukhan, P., and Doukhan, P.}
\newblock Mixing.
\newblock {\em Mixing: Properties and Examples\/} (1994), 15--23.

\bibitem{fan2024factor}
{\sc Fan, J., and Gu, Y.}
\newblock Factor augmented sparse throughput deep relu neural networks for high
  dimensional regression.
\newblock {\em Journal of the American Statistical Association 119}, 548
  (2024), 2680--2694.

\bibitem{fan2024noise}
{\sc Fan, J., Gu, Y., and Zhou, W.-X.}
\newblock How do noise tails impact on deep relu networks?
\newblock {\em The Annals of Statistics 52}, 4 (2024), 1845--1871.

\bibitem{fan2001variable}
{\sc Fan, J., and Li, R.}
\newblock Variable selection via nonconcave penalized likelihood and its oracle
  properties.
\newblock {\em Journal of the American statistical Association 96}, 456 (2001),
  1348--1360.

\bibitem{hang2016learning}
{\sc Hang, H., Feng, Y., Steinwart, I., and Suykens, J.~A.}
\newblock Learning theory estimates with observations from general stationary
  stochastic processes.
\newblock {\em Neural computation 28}, 12 (2016), 2853--2889.

\bibitem{hang2017bernstein}
{\sc Hang, H., and Steinwart, I.}
\newblock A bernstein-type inequality for some mixing processes and dynamical
  systems with an application to learning.
\newblock {\em The Annals of Statistics\/} (2017), 708--743.

\bibitem{jiao2023deep}
{\sc Jiao, Y., Shen, G., Lin, Y., and Huang, J.}
\newblock Deep nonparametric regression on approximate manifolds: Nonasymptotic
  error bounds with polynomial prefactors.
\newblock {\em The Annals of Statistics 51}, 2 (2023), 691--716.

\bibitem{juditsky2009nonparametric}
{\sc Juditsky, A.~B., Lepski, O., and Tsybakov, A.~B.}
\newblock Nonparametric estimation of composite functions.
\newblock {\em Annals of Statistics 37}, 3 (2009), 1360--1404.

\bibitem{kengne2025excess}
{\sc Kengne, W.}
\newblock Excess risk bound for deep learning under weak dependence.
\newblock {\em Mathematical Methods in the Applied Sciences\/} (2025).

\bibitem{kengne2025deep}
{\sc Kengne, W., and Wade, M.}
\newblock Deep learning from strongly mixing observations: Sparse-penalized
  regularization and minimax optimality.
\newblock {\em Journal of Complexity\/} (2025), 101978.

\bibitem{kengne2025robust}
{\sc Kengne, W., and Wade, M.}
\newblock Robust deep learning from weakly dependent data.
\newblock {\em Neural Networks 185\/} (2025), 107227.

\bibitem{kengne2025sparse}
{\sc Kengne, W., and Wade, M.}
\newblock Sparse-penalized deep neural networks estimator under weak
  dependence: W. kengne, m. wade.
\newblock {\em Metrika 88}, 4 (2025), 469--500.

\bibitem{kim2021fast}
{\sc Kim, Y., Ohn, I., and Kim, D.}
\newblock Fast convergence rates of deep neural networks for classification.
\newblock {\em Neural Networks 138\/} (2021), 179--197.

\bibitem{kohler2023rate}
{\sc Kohler, M., and Krzy{\.z}ak, A.}
\newblock On the rate of convergence of a deep recurrent neural network
  estimate in a regression problem with dependent data.
\newblock {\em Bernoulli 29}, 2 (2023), 1663--1685.

\bibitem{kurisu2025adaptive}
{\sc Kurisu, D., Fukami, R., and Koike, Y.}
\newblock Adaptive deep learning for nonlinear time series models.
\newblock {\em Bernoulli 31}, 1 (2025), 240--270.

\bibitem{maume2006exponential}
{\sc Maume-Deschamps, V.}
\newblock Exponential inequalities and functional estimations for weak
  dependent data: applications to dynamical systems.
\newblock {\em Stochastics and Dynamics 6}, 04 (2006), 535--560.

\bibitem{merlevede2009bernstein}
{\sc Merlevede, F., Peligrad, M., and Rio, E.}
\newblock Bernstein inequality and moderate deviations under strong mixing
  conditions.
\newblock In {\em High dimensional probability V: the Luminy volume}, vol.~5.
  Institute of Mathematical Statistics, 2009, pp.~273--293.

\bibitem{modha1996minimum}
{\sc Modha, D.~S., and Masry, E.}
\newblock Minimum complexity regression estimation with weakly dependent
  observations.
\newblock {\em IEEE Transactions on Information Theory 42}, 6 (1996),
  2133--2145.

\bibitem{ohn2019smooth}
{\sc Ohn, I., and Kim, Y.}
\newblock Smooth function approximation by deep neural networks with general
  activation functions.
\newblock {\em Entropy 21}, 7 (2019), 627.

\bibitem{ohn2022nonconvex}
{\sc Ohn, I., and Kim, Y.}
\newblock Nonconvex sparse regularization for deep neural networks and its
  optimality.
\newblock {\em Neural Computation 34}, 2 (2022), 476--517.

\bibitem{padilla2022quantile}
{\sc Padilla, O. H.~M., Tansey, W., and Chen, Y.}
\newblock Quantile regression with relu networks: Estimators and minimax rates.
\newblock {\em Journal of Machine Learning Research 23}, 247 (2022), 1--42.

\bibitem{samson2000concentration}
{\sc Samson, P.-M.}
\newblock Concentration of measure inequalities for markov chains and $phi
  $-mixing processes.
\newblock {\em The Annals of Probability 28}, 1 (2000), 416--461.

\bibitem{schmidt2020nonparametric}
{\sc Schmidt-Hieber, J.}
\newblock Nonparametric regression using deep neural networks with relu
  activation function.
\newblock {\em The Annals of Statistics 48}, 4 (2020), 1875--1897.

\bibitem{shen2021deep}
{\sc Shen, G., Jiao, Y., Lin, Y., Horowitz, J.~L., and Huang, J.}
\newblock Deep quantile regression: Mitigating the curse of dimensionality
  through composition.
\newblock {\em arXiv preprint arXiv:2107.04907\/} (2021).

\bibitem{shen2024nonparametric}
{\sc Shen, G., Jiao, Y., Lin, Y., Horowitz, J.~L., and Huang, J.}
\newblock Nonparametric estimation of non-crossing quantile regression process
  with deep requ neural networks.
\newblock {\em Journal of Machine Learning Research 25}, 88 (2024), 1--75.

\bibitem{zhang2010nearly}
{\sc Zhang, C.-H.}
\newblock Nearly unbiased variable selection under minimax concave penalty.
\newblock {\em The Annals of Statistics\/} (2010), 894--942.

\bibitem{zhang2010analysis}
{\sc Zhang, T.}
\newblock Analysis of multi-stage convex relaxation for sparse regularization.
\newblock {\em Journal of Machine Learning Research 11}, 3 (2010).

\bibitem{zhang2024classification}
{\sc Zhang, Z., Shi, L., and Zhou, D.-X.}
\newblock Classification with deep neural networks and logistic loss.
\newblock {\em Journal of Machine Learning Research 25}, 125 (2024), 1--117.

\end{thebibliography}
\end{document}